\pdfoutput=1

\documentclass[12pt,a4paper]{article}

\usepackage{amsfonts,amssymb,amsthm,amsmath,latexsym}
\usepackage{eqsection,indent}
\usepackage{subeqnarray, eepicemu, url,cite,bm}
\usepackage{multirow}  
\usepackage{tensor}  
\usepackage{graphicx,graphbox}
\usepackage{float}  

\date{April 25, 2022 \\[1mm] revised September 9, 2022}

\begin{document}

\title{\vspace*{-1cm}
       Multiple orthogonal polynomials, \\
       $d$-orthogonal polynomials, \\
       production matrices, \\
       and branched continued fractions
      }

\author{
     {\small Alan D.~Sokal}                  \\[2mm]
     {\small\it Department of Mathematics}   \\[-2mm]
     {\small\it University College London}   \\[-2mm]
     {\small\it Gower Street}                \\[-2mm]
     {\small\it London WC1E 6BT}             \\[-2mm]
     {\small\it UNITED KINGDOM}              \\[-2mm]
     {\small\tt sokal@math.ucl.ac.uk}        \\[-2mm]
     {\protect\makebox[5in]{\quad}}  
     \\[-2mm]
     {\small\it Department of Physics}       \\[-2mm]
     {\small\it New York University}         \\[-2mm]
     {\small\it 726 Broadway}                \\[-2mm]
     {\small\it New York, NY 10003}          \\[-2mm]
     {\small\it USA}                         \\[-2mm]
     {\small\tt sokal@nyu.edu}               \\[3mm]
}

\maketitle
\thispagestyle{empty}   

\begin{abstract}
I analyze an unexpected connection between
multiple orthogonal polynomials, $d$-orthogonal polynomials,
production matrices and branched continued fractions.
This work can be viewed as a partial extension of
Viennot's combinatorial theory of orthogonal polynomials
to the case where the production matrix is lower-Hessenberg
but is not necessarily tridiagonal.
\end{abstract}


\bigskip
\bigskip
\noindent
{\bf Mathematics Subject Classification (MSC 2010) codes:}
42C05 (Primary);
05A15, 05A19, 15A24, 15A99, 15B99, 30B70, 30E05, 33C45, 41A21 (Secondary).

\clearpage


\newtheorem{theorem}{Theorem}[section]
\newtheorem{proposition}[theorem]{Proposition}
\newtheorem{lemma}[theorem]{Lemma}
\newtheorem{corollary}[theorem]{Corollary}
\newtheorem{definition}[theorem]{Definition}
\newtheorem{conjecture}[theorem]{Conjecture}
\newtheorem{question}[theorem]{Question}
\newtheorem{problem}[theorem]{Problem}
\newtheorem{openproblem}[theorem]{Open Problem}
\newtheorem{example}[theorem]{Example}

\renewcommand{\theenumi}{\alph{enumi}}
\renewcommand{\labelenumi}{(\theenumi)}
\def\eop{\hbox{\kern1pt\vrule height6pt width4pt
depth1pt\kern1pt}\medskip}
\def\prf{\par\noindent{\bf Proof.\enspace}\rm}
\def\rmk{\par\medskip\noindent{\bf Remark\enspace}\rm}

\newcommand{\textbfit}[1]{\textbf{\textit{#1}}}

\newcommand{\bigdash}{%
\smallskip\begin{center} \rule{5cm}{0.1mm} \end{center}\smallskip}

\newcommand{\safepar}{ {\protect\hfill\protect\break\hspace*{5mm}} }

\newcommand{\be}{\begin{equation}}
\newcommand{\ee}{\end{equation}}
\newcommand{\<}{\langle}
\renewcommand{\>}{\rangle}
\newcommand{\widebar}{\overline}
\def\reff#1{(\protect\ref{#1})}
\def\spose#1{\hbox to 0pt{#1\hss}}
\def\ltapprox{\mathrel{\spose{\lower 3pt\hbox{$\mathchar"218$}}
    \raise 2.0pt\hbox{$\mathchar"13C$}}}
\def\gtapprox{\mathrel{\spose{\lower 3pt\hbox{$\mathchar"218$}}
    \raise 2.0pt\hbox{$\mathchar"13E$}}}
\def\textprime{${}^\prime$}
\def\proof{\par\medskip\noindent{\sc Proof.\ }}
\def\firstproof{\par\medskip\noindent{\sc First Proof.\ }}
\def\secondproof{\par\medskip\noindent{\sc Second Proof.\ }}
\def\alternateproof{\par\medskip\noindent{\sc Alternate Proof.\ }}
\def\algebraicproof{\par\medskip\noindent{\sc Algebraic Proof.\ }}
\def\combinatorialproof{\par\medskip\noindent{\sc Combinatorial Proof.\ }}
\def\proofof#1{\bigskip\noindent{\sc Proof of #1.\ }}
\def\firstproofof#1{\bigskip\noindent{\sc First Proof of #1.\ }}
\def\secondproofof#1{\bigskip\noindent{\sc Second Proof of #1.\ }}
\def\thirdproofof#1{\bigskip\noindent{\sc Third Proof of #1.\ }}
\def\algebraicproofof#1{\bigskip\noindent{\sc Algebraic Proof of #1.\ }}
\def\combinatorialproofof#1{\bigskip\noindent{\sc Combinatorial Proof of #1.\ }}
\def\sketchofproof{\par\medskip\noindent{\sc Sketch of proof.\ }}
\renewcommand{\qed}{ $\square$ \bigskip}
\newcommand{\myendremark}{ $\blacksquare$ \bigskip}
\def\half{ {1 \over 2} }
\def\third{ {1 \over 3} }
\def\twothird{ {2 \over 3} }
\def\smfrac#1#2{{\textstyle{#1\over #2}}}
\def\smhalf{ {\smfrac{1}{2}} }
\newcommand{\real}{\mathop{\rm Re}\nolimits}
\renewcommand{\Re}{\mathop{\rm Re}\nolimits}
\newcommand{\imag}{\mathop{\rm Im}\nolimits}
\renewcommand{\Im}{\mathop{\rm Im}\nolimits}
\newcommand{\sgn}{\mathop{\rm sgn}\nolimits}
\newcommand{\tr}{\mathop{\rm tr}\nolimits}
\newcommand{\supp}{\mathop{\rm supp}\nolimits}
\newcommand{\disc}{\mathop{\rm disc}\nolimits}
\newcommand{\diag}{\mathop{\rm diag}\nolimits}
\newcommand{\tridiag}{\mathop{\rm tridiag}\nolimits}
\newcommand{\AZ}{\mathop{\rm AZ}\nolimits}
\newcommand{\NC}{\mathop{\rm NC}\nolimits}
\newcommand{\PF}{{\rm PF}}
\newcommand{\rk}{\mathop{\rm rk}\nolimits}
\newcommand{\perm}{\mathop{\rm perm}\nolimits}
\def\hboxscript#1{ {\hbox{\scriptsize\em #1}} }
\renewcommand{\emptyset}{\varnothing}
\newcommand{\eqdef}{\stackrel{\rm def}{=}}

\newcommand{\restrict}{\upharpoonright}

\newcommand{\compinv}{{\langle -1 \rangle}}   

\newcommand{\divides}{\mid}
\newcommand{\notdivides}{\nmid}

\newcommand{\scra}{{\mathcal{A}}}
\newcommand{\scrb}{{\mathcal{B}}}
\newcommand{\scrc}{{\mathcal{C}}}
\newcommand{\scrd}{{\mathcal{D}}}
\newcommand{\scrdtilde}{{\widetilde{\mathcal{D}}}}
\newcommand{\scre}{{\mathcal{E}}}
\newcommand{\scrf}{{\mathcal{F}}}
\newcommand{\scrg}{{\mathcal{G}}}
\newcommand{\scrh}{{\mathcal{H}}}
\newcommand{\scri}{{\mathcal{I}}}
\newcommand{\scrj}{{\mathcal{J}}}
\newcommand{\scrk}{{\mathcal{K}}}
\newcommand{\scrl}{{\mathcal{L}}}
\newcommand{\scrm}{{\mathcal{M}}}
\newcommand{\scrn}{{\mathcal{N}}}
\newcommand{\scro}{{\mathcal{O}}}
\newcommand\scroo{
  \mathchoice
    {{\scriptstyle\mathcal{O}}}
    {{\scriptstyle\mathcal{O}}}
    {{\scriptscriptstyle\mathcal{O}}}
    {\scalebox{0.6}{$\scriptscriptstyle\mathcal{O}$}}
  }
\newcommand{\scrp}{{\mathcal{P}}}
\newcommand{\scrq}{{\mathcal{Q}}}
\newcommand{\scrr}{{\mathcal{R}}}
\newcommand{\scrs}{{\mathcal{S}}}
\newcommand{\scrt}{{\mathcal{T}}}
\newcommand{\scrv}{{\mathcal{V}}}
\newcommand{\scrw}{{\mathcal{W}}}
\newcommand{\scrz}{{\mathcal{Z}}}
\newcommand{\SP}{{\mathcal{SP}}}
\newcommand{\ST}{{\mathcal{ST}}}

\newcommand{\bfa}{{\mathbf{a}}}
\newcommand{\bfb}{{\mathbf{b}}}
\newcommand{\bfc}{{\mathbf{c}}}
\newcommand{\bfd}{{\mathbf{d}}}
\newcommand{\bfe}{{\mathbf{e}}}
\newcommand{\bfh}{{\mathbf{h}}}
\newcommand{\bfj}{{\mathbf{j}}}
\newcommand{\bfi}{{\mathbf{i}}}
\newcommand{\bfk}{{\mathbf{k}}}
\newcommand{\bfl}{{\mathbf{l}}}
\newcommand{\bfL}{{\mathbf{L}}}
\newcommand{\bfm}{{\mathbf{m}}}
\newcommand{\bfn}{{\mathbf{n}}}
\newcommand{\bfp}{{\mathbf{p}}}
\newcommand{\bfP}{{\mathbf{P}}}
\newcommand{\bfr}{{\mathbf{r}}}
\newcommand{\bfs}{{\mathbf{s}}}
\newcommand{\bfu}{{\mathbf{u}}}
\newcommand{\bfv}{{\mathbf{v}}}
\newcommand{\bfw}{{\mathbf{w}}}
\newcommand{\bfx}{{\mathbf{x}}}
\newcommand{\bfX}{{\mathbf{X}}}
\newcommand{\bfy}{{\mathbf{y}}}
\newcommand{\bfz}{{\mathbf{z}}}
\renewcommand{\k}{{\mathbf{k}}}
\newcommand{\n}{{\mathbf{n}}}
\newcommand{\vv}{{\mathbf{v}}}
\newcommand{\bv}{{\mathbf{v}}}
\newcommand{\w}{{\mathbf{w}}}
\newcommand{\x}{{\mathbf{x}}}
\newcommand{\y}{{\mathbf{y}}}
\newcommand{\cc}{{\mathbf{c}}}
\newcommand{\zero}{{\mathbf{0}}}
\newcommand{\one}{{\mathbf{1}}}
\newcommand{\bmm}{{\mathbf{m}}}

\newcommand{\ahat}{{\widehat{a}}}
\newcommand{\Zhat}{{\widehat{Z}}}

\newcommand{\C}{{\mathbb C}}
\newcommand{\D}{{\mathbb D}}
\newcommand{\Z}{{\mathbb Z}}
\newcommand{\N}{{\mathbb N}}
\newcommand{\Q}{{\mathbb Q}}
\newcommand{\PP}{{\mathbb P}}
\newcommand{\R}{{\mathbb R}}
\newcommand{\RR}{{\mathbb R}}
\newcommand{\E}{{\mathbb E}}

\newcommand{\Sym}{{\mathfrak{S}}}
\newcommand{\SymB}{{\mathfrak{B}}}
\newcommand{\Alt}{{\mathrm{Alt}}}

\newcommand{\germanA}{{\mathfrak{A}}}
\newcommand{\germanB}{{\mathfrak{B}}}
\newcommand{\germanQ}{{\mathfrak{Q}}}
\newcommand{\germanh}{{\mathfrak{h}}}

\newcommand{\myle}{\preceq}
\newcommand{\myge}{\succeq}
\newcommand{\mygt}{\succ}

\newcommand{\B}{{\sf B}}
\newcommand{\OB}{B^{\rm ord}}
\newcommand{\OS}{{\sf OS}}
\newcommand{\OO}{{\sf O}}
\newcommand{\OSP}{{\sf OSP}}
\newcommand{\Eu}{{\sf Eu}}
\newcommand{\ERR}{{\sf ERR}}
\newcommand{\sfB}{{\sf B}}
\newcommand{\sfD}{{\sf D}}
\newcommand{\sfE}{{\sf E}}
\newcommand{\sfG}{{\sf G}}
\newcommand{\sfJ}{{\sf J}}
\newcommand{\sfL}{{\sf L}}
\newcommand{\sfLhat}{{\widehat{{\sf L}}}}
\newcommand{\sfLtilde}{{\widetilde{{\sf L}}}}
\newcommand{\sfP}{{\sf P}}
\newcommand{\sfQ}{{\sf Q}}
\newcommand{\sfS}{{\sf S}}
\newcommand{\sfT}{{\sf T}}
\newcommand{\sfW}{{\sf W}}
\newcommand{\sfMV}{{\sf MV}}
\newcommand{\AMV}{{\sf AMV}}
\newcommand{\BM}{{\sf BM}}
\newcommand{\emIB}{B^{\rm irr}}
\newcommand{\emIP}{P^{\rm irr}}
\newcommand{\emOB}{B^{\rm ord}}
\newcommand{\emCB}{B^{\rm cyc}}
\newcommand{\emSC}{P^{\rm cyc}}

\newcommand{\lev}{{\rm lev}}
\newcommand{\stat}{{\rm stat}}
\newcommand{\cyc}{{\rm cyc}}
\newcommand{\mysteryone}{{\rm mys1}}
\newcommand{\mysterytwo}{{\rm mys2}}
\newcommand{\Asc}{{\rm Asc}}
\newcommand{\asc}{{\rm asc}}
\newcommand{\Des}{{\rm Des}}
\newcommand{\des}{{\rm des}}
\newcommand{\Exc}{{\rm Exc}}
\newcommand{\exc}{{\rm exc}}
\newcommand{\Wex}{{\rm Wex}}
\newcommand{\wex}{{\rm wex}}
\newcommand{\Fix}{{\rm Fix}}
\newcommand{\fix}{{\rm fix}}
\newcommand{\lrmax}{{\rm lrmax}}
\newcommand{\rlmax}{{\rm rlmax}}
\newcommand{\Rec}{{\rm Rec}}
\newcommand{\rec}{{\rm rec}}
\newcommand{\Arec}{{\rm Arec}}
\newcommand{\arec}{{\rm arec}}
\newcommand{\ERec}{{\rm ERec}}
\newcommand{\erec}{{\rm erec}}
\newcommand{\EArec}{{\rm EArec}}
\newcommand{\earec}{{\rm earec}}
\newcommand{\recarec}{{\rm recarec}}
\newcommand{\nonrec}{{\rm nonrec}}
\newcommand{\Cpeak}{{\rm Cpeak}}
\newcommand{\cpeak}{{\rm cpeak}}
\newcommand{\Cval}{{\rm Cval}}
\newcommand{\cval}{{\rm cval}}
\newcommand{\Cdasc}{{\rm Cdasc}}
\newcommand{\cdasc}{{\rm cdasc}}
\newcommand{\Cddes}{{\rm Cddes}}
\newcommand{\cddes}{{\rm cddes}}
\newcommand{\cdrise}{{\rm cdrise}}
\newcommand{\cdfall}{{\rm cdfall}}
\newcommand{\Peak}{{\rm Peak}}
\newcommand{\peak}{{\rm peak}}
\newcommand{\Val}{{\rm Val}}
\newcommand{\val}{{\rm val}}
\newcommand{\Dasc}{{\rm Dasc}}
\newcommand{\dasc}{{\rm dasc}}
\newcommand{\Ddes}{{\rm Ddes}}
\newcommand{\ddes}{{\rm ddes}}
\newcommand{\inv}{{\rm inv}}
\newcommand{\maj}{{\rm maj}}
\newcommand{\rs}{{\rm rs}}
\newcommand{\cross}{{\rm cr}}
\newcommand{\crosshat}{{\widehat{\rm cr}}}
\newcommand{\nest}{{\rm ne}}
\newcommand{\rodd}{{\rm rodd}}
\newcommand{\reven}{{\rm reven}}
\newcommand{\lodd}{{\rm lodd}}
\newcommand{\leven}{{\rm leven}}
\newcommand{\sg}{{\rm sg}}
\newcommand{\bl}{{\rm bl}}
\newcommand{\tran}{{\rm tr}}
\newcommand{\area}{{\rm area}}
\newcommand{\ret}{{\rm ret}}
\newcommand{\peaks}{{\rm peaks}}
\newcommand{\hl}{{\rm hl}}
\newcommand{\sll}{{\rm sl}}
\newcommand{\negg}{{\rm neg}}
\newcommand{\imp}{{\rm imp}}
\newcommand{\osg}{{\rm osg}}
\newcommand{\ons}{{\rm ons}}
\newcommand{\isg}{{\rm isg}}
\newcommand{\ins}{{\rm ins}}
\newcommand{\LL}{{\rm LL}}
\newcommand{\height}{{\rm ht}}
\newcommand{\as}{{\rm as}}

\newcommand{\ba}{{\bm{a}}}
\newcommand{\bahat}{{\widehat{\bm{a}}}}
\newcommand{\sfa}{{{\sf a}}}
\newcommand{\bb}{{\bm{b}}}
\newcommand{\bc}{{\bm{c}}}
\newcommand{\bchat}{{\widehat{\bm{c}}}}
\newcommand{\bd}{{\bm{d}}}
\newcommand{\bee}{{\bm{e}}}
\newcommand{\beh}{{\bm{eh}}}
\newcommand{\bff}{{\bm{f}}}
\newcommand{\bg}{{\bm{g}}}
\newcommand{\bh}{{\bm{h}}}
\newcommand{\bll}{{\bm{\ell}}}
\newcommand{\bp}{{\bm{p}}}
\newcommand{\br}{{\bm{r}}}
\newcommand{\bs}{{\bm{s}}}
\newcommand{\bu}{{\bm{u}}}
\newcommand{\bw}{{\bm{w}}}
\newcommand{\bx}{{\bm{x}}}
\newcommand{\by}{{\bm{y}}}
\newcommand{\bz}{{\bm{z}}}
\newcommand{\bA}{{\bm{A}}}
\newcommand{\bB}{{\bm{B}}}
\newcommand{\bC}{{\bm{C}}}
\newcommand{\bE}{{\bm{E}}}
\newcommand{\bF}{{\bm{F}}}
\newcommand{\bG}{{\bm{G}}}
\newcommand{\bH}{{\bm{H}}}
\newcommand{\bI}{{\bm{I}}}
\newcommand{\bJ}{{\bm{J}}}
\newcommand{\bL}{{\bm{L}}}
\newcommand{\bLhat}{{\widehat{\bm{L}}}}
\newcommand{\bM}{{\bm{M}}}
\newcommand{\bN}{{\bm{N}}}
\newcommand{\bP}{{\bm{P}}}
\newcommand{\bQ}{{\bm{Q}}}
\newcommand{\bR}{{\bm{R}}}
\newcommand{\bS}{{\bm{S}}}
\newcommand{\bT}{{\bm{T}}}
\newcommand{\bW}{{\bm{W}}}
\newcommand{\bX}{{\bm{X}}}
\newcommand{\bY}{{\bm{Y}}}
\newcommand{\bIB}{{\bm{B}^{\rm irr}}}
\newcommand{\bOB}{{\bm{B}^{\rm ord}}}
\newcommand{\bOS}{{\bm{OS}}}
\newcommand{\bERR}{{\bm{ERR}}}
\newcommand{\bSP}{{\bm{SP}}}
\newcommand{\bMV}{{\bm{MV}}}
\newcommand{\bBM}{{\bm{BM}}}
\newcommand{\balpha}{{\bm{\alpha}}}
\newcommand{\balphapre}{{\bm{\alpha}^{\rm pre}}}
\newcommand{\bbeta}{{\bm{\beta}}}
\newcommand{\bgamma}{{\bm{\gamma}}}
\newcommand{\bGamma}{{\bm{\Gamma}}}
\newcommand{\bdelta}{{\bm{\delta}}}
\newcommand{\bkappa}{{\bm{\kappa}}}
\newcommand{\bmu}{{\bm{\mu}}}
\newcommand{\bomega}{{\bm{\omega}}}
\newcommand{\bsigma}{{\bm{\sigma}}}
\newcommand{\btau}{{\bm{\tau}}}
\newcommand{\bphi}{{\bm{\phi}}}
\newcommand{\bphihat}{{\skew{3}\widehat{\vphantom{t}\protect\smash{\bm{\phi}}}}}
\newcommand{\bpsi}{{\bm{\psi}}}
\newcommand{\bxi}{{\bm{\xi}}}
\newcommand{\bzeta}{{\bm{\zeta}}}
\newcommand{\bone}{{\bm{1}}}
\newcommand{\bzero}{{\bm{0}}}

\newcommand{\Cbar}{{\overline{C}}}
\newcommand{\Dbar}{{\overline{D}}}
\newcommand{\dbar}{{\overline{d}}}
\def\Btilde{{\widetilde{B}}}
\def\Ctilde{{\widetilde{C}}}
\def\Ftilde{{\widetilde{F}}}
\def\Gtilde{{\widetilde{G}}}
\def\Htilde{{\widetilde{H}}}
\def\Lhat{{\widehat{L}}}
\def\Ltilde{{\widetilde{L}}}
\def\Ptilde{{\widetilde{P}}}
\def\Phat{{\widehat{P}}}
\def\bfPhat{{\widehat{\bfP}}}
\def\ptilde{{\widetilde{p}}}
\def\Chat{{\widehat{C}}}
\def\ctilde{{\widetilde{c}}}
\def\zbar{{\overline{Z}}}
\def\pitilde{{\widetilde{\pi}}}
\def\omegahat{{\widehat{\omega}}}

\newcommand{\sech}{{\rm sech}}

%
%
\newcommand{\sn}{{\rm sn}}
\newcommand{\cn}{{\rm cn}}
\newcommand{\dn}{{\rm dn}}
\newcommand{\sm}{{\rm sm}}
\newcommand{\cm}{{\rm cm}}

%
%
\newcommand{\zfz}{ {{}_0 \! F_0} }
\newcommand{\zfo}{ {{}_0  F_1} }
\newcommand{\ofz}{ {{}_1 \! F_0} }
\newcommand{\ofo}{ {{}_1 \! F_1} }
\newcommand{\oft}{ {{}_1 \! F_2} }

%
%
\newcommand{\FHyper}[2]{ {\tensor[_{#1 \!}]{F}{_{#2}}\!} }
\newcommand{\FHYPER}[5]{ {\FHyper{#1}{#2} \!\biggl(
   \!\!\begin{array}{c} #3 \\[1mm] #4 \end{array}\! \bigg|\, #5 \! \biggr)} }
\newcommand{\tfo}{ {\FHyper{2}{1}} }
\newcommand{\tfz}{ {\FHyper{2}{0}} }
\newcommand{\threefz}{ {\FHyper{3}{0}} }
\newcommand{\FHYPERbottomzero}[3]{ {\FHyper{#1}{0} \hspace*{-0mm}\biggl(
   \!\!\begin{array}{c} #2 \\[1mm] \hbox{---} \end{array}\! \bigg|\, #3 \! \biggr)} }
\newcommand{\FHYPERtopzero}[3]{ {\FHyper{0}{#1} \hspace*{-0mm}\biggl(
   \!\!\begin{array}{c} \hbox{---} \\[1mm] #2 \end{array}\! \bigg|\, #3 \! \biggr)} }

\newcommand{\phiHyper}[2]{ {\tensor[_{#1}]{\phi}{_{#2}}} }
\newcommand{\psiHyper}[2]{ {\tensor[_{#1}]{\psi}{_{#2}}} }
\newcommand{\PhiHyper}[2]{ {\tensor[_{#1}]{\Phi}{_{#2}}} }
\newcommand{\PsiHyper}[2]{ {\tensor[_{#1}]{\Psi}{_{#2}}} }
\newcommand{\phiHYPER}[6]{ {\phiHyper{#1}{#2} \!\left(
   \!\!\begin{array}{c} #3 \\ #4 \end{array}\! ;\, #5, \, #6 \! \right)\!} }
\newcommand{\psiHYPER}[6]{ {\psiHyper{#1}{#2} \!\left(
   \!\!\begin{array}{c} #3 \\ #4 \end{array}\! ;\, #5, \, #6 \! \right)} }
\newcommand{\PhiHYPER}[5]{ {\PhiHyper{#1}{#2} \!\left(
   \!\!\begin{array}{c} #3 \\ #4 \end{array}\! ;\, #5 \! \right)\!} }
\newcommand{\PsiHYPER}[5]{ {\PsiHyper{#1}{#2} \!\left(
   \!\!\begin{array}{c} #3 \\ #4 \end{array}\! ;\, #5 \! \right)\!} }
\newcommand{\zerophizero}{ {\phiHyper{0}{0}} }
\newcommand{\ophizero}{ {\phiHyper{1}{0}} }
\newcommand{\zphio}{ {\phiHyper{0}{1}} }
\newcommand{\ophio}{ {\phiHyper{1}{1}} }
\newcommand{\tphio}{ {\phiHyper{2}{1}} }
\newcommand{\tphiz}{ {\phiHyper{2}{0}} }
\newcommand{\tPhio}{ {\PhiHyper{2}{1}} }
\newcommand{\opsio}{ {\psiHyper{1}{1}} }

%
%
\newcommand{\stirlingsubset}[2]{\genfrac{\{}{\}}{0pt}{}{#1}{#2}}
\newcommand{\stirlingcycle}[2]{\genfrac{[}{]}{0pt}{}{#1}{#2}}
\newcommand{\assocstirlingsubset}[3]{{\genfrac{\{}{\}}{0pt}{}{#1}{#2}}_{\! \ge #3}}
\newcommand{\genstirlingsubset}[4]{{\genfrac{\{}{\}}{0pt}{}{#1}{#2}}_{\! #3,#4}}
\newcommand{\irredstirlingsubset}[2]{{\genfrac{\{}{\}}{0pt}{}{#1}{#2}}^{\!\rm irr}}
\newcommand{\euler}[2]{\genfrac{\langle}{\rangle}{0pt}{}{#1}{#2}}
\newcommand{\eulergen}[3]{{\genfrac{\langle}{\rangle}{0pt}{}{#1}{#2}}_{\! #3}}
\newcommand{\eulersecond}[2]{\left\langle\!\! \euler{#1}{#2} \!\!\right\rangle}
\newcommand{\eulersecondgen}[3]{{\left\langle\!\! \euler{#1}{#2} \!\!\right\rangle}_{\! #3}}
\newcommand{\binomvert}[2]{\genfrac{\vert}{\vert}{0pt}{}{#1}{#2}}
\newcommand{\binomsquare}[2]{\genfrac{[}{]}{0pt}{}{#1}{#2}}
\newcommand{\doublebinom}[2]{\left(\!\! \binom{#1}{#2} \!\!\right)}


\newenvironment{sarray}{
             \textfont0=\scriptfont0
             \scriptfont0=\scriptscriptfont0
             \textfont1=\scriptfont1
             \scriptfont1=\scriptscriptfont1
             \textfont2=\scriptfont2
             \scriptfont2=\scriptscriptfont2
             \textfont3=\scriptfont3
             \scriptfont3=\scriptscriptfont3
           \renewcommand{\arraystretch}{0.7}
           \begin{array}{l}}{\end{array}}

\newenvironment{scarray}{
             \textfont0=\scriptfont0
             \scriptfont0=\scriptscriptfont0
             \textfont1=\scriptfont1
             \scriptfont1=\scriptscriptfont1
             \textfont2=\scriptfont2
             \scriptfont2=\scriptscriptfont2
             \textfont3=\scriptfont3
             \scriptfont3=\scriptscriptfont3
           \renewcommand{\arraystretch}{0.7}
           \begin{array}{c}}{\end{array}}


\newcommand*\circled[1]{\tikz[baseline=(char.base)]{
  \node[shape=circle,draw,inner sep=1pt] (char) {#1};}}
\newcommand{\ostar}{{\circledast}}
\newcommand{\ostarN}{{\,\circledast_{\vphantom{\dot{N}}N}\,}}
\newcommand{\ostarPsi}{{\,\circledast_{\vphantom{\dot{\Psi}}\Psi}\,}}
\newcommand{\starN}{{\,\ast_{\vphantom{\dot{N}}N}\,}}
\newcommand{\starpsi}{{\,\ast_{\vphantom{\dot{\bpsi}}\!\bpsi}\,}}
\newcommand{\starone}{{\,\ast_{\vphantom{\dot{1}}1}\,}}
\newcommand{\startwo}{{\,\ast_{\vphantom{\dot{2}}2}\,}}
\newcommand{\starinfty}{{\,\ast_{\vphantom{\dot{\infty}}\infty}\,}}
\newcommand{\starT}{{\,\ast_{\vphantom{\dot{T}}T}\,}}

\newcommand*{\Scale}[2][4]{\scalebox{#1}{$#2$}}

\newcommand*{\Scaletext}[2][4]{\scalebox{#1}{#2}} 

\tableofcontents

\clearpage

\section{Introduction}

The goal of this paper is to point out, and then analyze in detail,
an unexpected connection between multiple orthogonal polynomials
and $d$-orthogonal polynomials on the one hand,
and production matrices and branched continued fractions on the other
--- objects that arose over the past few decades
in the special-functions and enumerative-combinatorics communities,
respectively.
It is appropriate to begin, therefore,
by explaining briefly each of these four concepts.

Multiple orthogonal polynomials
\cite{Aptekarev_98,Martinez-Finkelshtein_16,VanAssche_21}
\cite[Chapter~23]{Ismail_05}
are a generalization of conventional orthogonal polynomials
\cite{Szego_75,Chihara_78,Ismail_05}
in which the polynomials satisfy orthogonality relations
with respect to several measures $\mu_1,\ldots,\mu_r$
rather than just one.
Multiple orthogonal polynomials first arose in Hermite--Pad\'e approximation
\cite[Chapter~4]{Nikishin_91} \cite{VanAssche_06};
they have applications to number theory
\cite{Sorokin_94,Sorokin_02,VanAssche_10},
random matrices \cite{Kuijlaars_10a,Kuijlaars_10b,Kuijlaars_14}
and nonintersecting random paths
\cite{Kuijlaars_09,Kuijlaars_10a,Kuijlaars_10b},
among other fields.
See Section~\ref{subsec.MOP} for a brief summary.

Closely related to multiple orthogonal polynomials
are the so-called $d$-orthogonal polynomials.
A sequence $(P_n(x))_{n \ge 0}$ of monic polynomials
is said to be $d$-orthogonal \cite{VanIseghem_87,Maroni_89}
with respect to a sequence $\Gamma_0,\ldots,\Gamma_{d-1}$ of linear forms
in case ${\Gamma_k(x^\ell \, P_n(x)) = 0}$ whenever $n > d\ell+k$.
(To avoid trivialities, it is also usually required that
 ${\Gamma_k(x^\ell \, P_n(x)) \neq 0}$ when $n = d\ell+k$.)
For $d=1$ this reduces to the ordinary concept of orthogonality.
It turns out that the sequence of multiple orthogonal polynomials of type~II
taken along the so-called stepline
is $d$-orthogonal (for $d=r$) with respect to
the linear forms associated to the measures $\mu_1,\ldots,\mu_r$.

Production matrices \cite{Deutsch_05,Deutsch_09}
have become in recent years an important tool in enumerative combinatorics
(see Section~\ref{subsec.prodmat_intro} for a brief summary).
In the special case of a tridiagonal production matrix,
this construction goes back to Stieltjes' \cite{Stieltjes_1889,Stieltjes_1894}
work on continued fractions:
the production matrix of a classical S-fraction or J-fraction is tridiagonal.
Moreover, the classical J-fraction and tridiagonal production matrix
associated to the moment sequence of a measure $\mu$
are closely related to the sequence of orthogonal polynomials
associated to $\mu$.
This connection was comprehensively investigated
by Viennot \cite{Viennot_83,Viennot_85} in the early 1980s,
who developed a general combinatorial theory of orthogonal polynomials,
building on Flajolet's \cite{Flajolet_80} combinatorial theory of
continued fractions.
Our work here can be viewed as a partial extension of Viennot's theory
to the case where the production matrix is lower-Hessenberg
(i.e.\ vanishes above the first superdiagonal)
but is not necessarily tridiagonal.
Indeed, this extension was already begun by Viennot himself
\cite[sections~III.5 and V.6]{Viennot_83}.

Various types of branched continued fractions have been introduced
in the analysis literature
\cite{DeBruin_78,Bodnar_98,Bodnar_18}
\cite[pp.~274--280, 285]{Lorentzen_92} \cite[p.~28]{Cuyt_08},
but we are not concerned here with these.
Rather, we are concerned with the branched continued fractions
that have been introduced by combinatorialists
and whose Taylor coefficients are the generating polynomials
for selected types of lattice paths,
generalizing the work of Flajolet \cite{Flajolet_80}
on classical continued fractions.
This investigation was also initiated
by Viennot \cite[section~V.6]{Viennot_83},
who briefly considered the branched continued fractions
({\em fractions multicontinu\'ees}\/)
generated by \L{}ukasiewicz paths.
This work was carried forward in the Ph.D.~theses of
Roblet \cite{Roblet_94} and Varvak \cite{Varvak_04}.
Further applications were made by
Gouyou-Beauchamps \cite{Gouyou-Beauchamps_98} and Drake \cite{Drake_11}.
Subsequently, Albenque and Bouttier \cite{Albenque_12}
introduced the branched continued fractions generated by $m$-Dyck paths
and proved many interesting results about them.
Most recently, P\'etr\'eolle, Sokal and Zhu \cite{latpath_SRTR}
carried out a comprehensive analysis of the branched continued fractions
associated to $m$-Dyck, $m$-Schr\"oder and \L{}ukasiewicz paths,
with emphasis on questions related to total positivity;
see also \cite{latpath_lah,latpath_laguerre} for further applications.
For these branched continued fractions, the production matrix
is lower-Hessenberg but not (except in the classical cases) tridiagonal.
See Section~\ref{subsec.BCF} for a  brief summary.

Finally, let us mention the remarkable Ph.D.~thesis of Drake \cite{Drake_06},
who initiated the combinatorial theory of multiple orthogonal polynomials
and who foresaw the link with branched continued fractions
\cite[p.~1]{Drake_06}.
Our work here can be viewed as an extension,
and to some extent a completion, of his.

Throughout this paper, we fix a commutative ring
(with identity element $1 \neq 0$) $R$:
we~will use sequences and matrices with entries in $R$,
and polynomials and formal power series with coefficients in $R$.
The analyst reader should feel free to imagine,
without too much loss of generality, that $R = \R$.
However, in applications of this formalism there will often be parameters,
and we will usually prefer to treat these parameters as
algebraic indeterminates~$\bxi$;  then $R$ will be either
the ring $\R[\bxi]$ of polynomials in these indeterminates
or the field $\R(\bxi)$ of rational functions in these indeterminates.

In particular, to any (positive or signed) measure $\mu$ on $\R$
that has finite moments of all orders, there is canonically associated
a linear functional $\scrl$ on the polynomial ring $\R[x]$,
defined by $\scrl(x^n) = \int\! x^n \, d\mu(x)$.
It is well known \cite{Chihara_78} that the theory of orthogonal polynomials
(or at least the simplest parts of it)
can be expressed entirely in terms of this linear functional
--- or equivalently, in terms of the sequence of moments $a_n = \scrl(x^n)$ ---
without reference to the measure $\mu$.
We shall adopt this approach here, and also replace $\R$
by an arbitrary commutative ring $R$.

The plan of this paper is as follows:
In Section~\ref{sec.prelim} we collect some basic definitions and facts
concerning multiple orthogonal polynomials, production matrices
and branched continued fractions.
In Section~\ref{sec.prodmat} we use the theory of production matrices
to demonstrate some very simple relations between
sequences of monic polynomials, the linear recurrences they satisfy,
and their dual sequences of linear functionals.
We also generalize Viennot's \cite{Viennot_83} formula
for the expectation of products of orthogonal polynomials.
In Section~\ref{sec.orthogonal}
we analyze sequences of monic polynomials that are orthogonal
to a sequence of linear functionals.
In Section~\ref{sec.prodmat_OOP} we show how this theory
applies to ordinary orthogonal polynomials,
and in Section~\ref{sec.prodmat_MOP}
we apply it to multiple orthogonal polynomials.
Finally, in Section~\ref{sec.examples} we examine some concrete examples.
In the Appendix we prove a basic result concerning
$LU$ factorization for matrices over a commutative ring.

\section{Preliminaries}   \label{sec.prelim}

In this section we provide a brief introduction to
multiple orthogonal polynomials
\cite{Aptekarev_98,Martinez-Finkelshtein_16,VanAssche_21}
\cite[Chapter~23]{Ismail_05},
production matrices \cite{Deutsch_05,Deutsch_09,latpath_SRTR,Sokal_totalpos},
and branched continued fractions \cite{latpath_SRTR}.
The reader familiar with one or more of these topics
can skim those parts quickly, with the main aim of fixing the notation.

\subsection{Multiple orthogonal polynomials}  \label{subsec.MOP}

We begin by giving a brief introduction to the theory of
multiple orthogonal polynomials,
following \cite[Chapter~23]{Ismail_05}
but making a few comments about ``algebraizing'' the theory
to allow coefficients in an arbitrary commutative ring $R$.
We limit attention to the multiple orthogonal polynomials of type~II,
since these are the only ones that will arise in the remainder of the paper.
(I leave it to others to investigate whether there is any analogue
of the connections discussed here for the
multiple orthogonal polynomials of type~I.)

Fix an integer $r \ge 1$,
and let $\mu_1,\ldots,\mu_r$ be positive measures on the real line
with finite moments of all orders.
We use multi-indices $\bfn = (n_1,\ldots,n_r) \in \N^r$
and write $|\bfn| = n_1 + \ldots + n_r$.
The \textbfit{multiple orthogonal polynomial of type~II}
for the multi-index $\bfn$ is the degree-$|\bfn|$ monic polynomial
$P_\bfn(x) = x^{|\bfn|} + \ldots\,$ satisfying
the orthogonality relations
\be
   \int\! x^k \, P_\bfn(x) \, d\mu_j(x)  \;=\; 0
   \quad
   \hbox{for all $1 \le j \le r$ and $0 \le k \le n_j - 1$}
   \;,
 \label{def.orthogonality}
\ee
whenever such a polynomial exists and is unique.
The equations \reff{def.orthogonality} give a system of $|\bfn|$
linear equations for the $|\bfn|$ non-leading coefficients of
the polynomial $P_\bfn(x)$;
the multi-index $\bfn$ is said to be {\em normal}\/
whenever the solution exists and is unique.
Note that the coefficient matrix of this system is the
transpose of the $|\bfn| \times |\bfn|$ matrix
\be
   \scrm_\bfn
   \;=\;
   \begin{pmatrix}
       H^{(1)}_{|\bfn|,n_1} &
       H^{(2)}_{|\bfn|,n_2} &
       \cdots &
       H^{(r)}_{|\bfn|,n_r}
   \end{pmatrix}
\ee
where $H^{(j)}_{M,N}$ is the $M \times N$ Hankel matrix
of the moments of $\mu_j$:  that is,
$H^{(j)}_{M,N} = (m^{(j)}_{r+s})_{0 \le r \le M-1,\, 0 \le s \le N-1}$
where
\be
   m^{(j)}_k  \;=\;  \int\! x^k \, d\mu_j(x)
   \;.
\ee
Therefore, the multi-index $\bfn$ is normal
if and only if $\det \scrm_\bfn \neq 0$.

A system of measures $\mu_1,\ldots,\mu_r$ is said to be {\em perfect}\/
in case all $\bfn \in \N^r$ are normal.
Several general sufficient conditions for a system to be perfect are known
(Angelesco systems, AT systems, Nikishin systems, \ldots):
see \cite[Chapter~23]{Ismail_05} \cite{VanAssche_21}.
We shall henceforth restrict attention to perfect systems.

The orthogonality conditions \reff{def.orthogonality}
can be trivially re-expressed in terms of
the linear forms $\scrl^{(1)},\ldots,\scrl^{(r)}$
associated to the measures $\mu_1,\ldots,\mu_r$,
which are defined by $\scrl^{(j)}(x^n) = \int\! x^n \, d\mu_j(x)$:
it suffices to replace $m^{(j)}_k$ by $\scrl^{(j)}(x^k)$.
Moreover, from this point of view,
the measures $\mu_j$ need not be positive measures;
indeed, the linear forms $\scrl^{(j)}$
need not come from (signed) measures at all.
Provided that one can show, one way or another,
that all $\bfn \in \N^r$ are normal,
the multiple orthogonal polynomials are well-defined.

Having done this, we can go farther and ``algebraize'' the theory
by considering polynomials with coefficients in an arbitrary commutative ring
(with identity element $1 \neq 0$) $R$, rather than just $R = \R$.
We fix linear forms $\scrl^{(1)},\ldots,\scrl^{(r)}$
on the polynomial ring $R[x]$,
and define ``moments'' $m^{(j)}_k = \scrl^{(j)}(x^k)$;
then the multi-index $\bfn$ is normal
if and only if $\det \scrm_\bfn$ is an invertible element of the ring $R$.

Let us now make a simple but important observation.
Fix a multi-index $\bfn = (n_1,\ldots,n_r) \in \N^r$,
and suppose that the polynomial $P_\bfn(x)$ satisfies the
orthogonality relations \reff{def.orthogonality}
with respect to some family of (not necessarily positive)
measures $\bmu = (\mu_1,\ldots,\mu_r)$.
Then $P_\bfn(x)$ also satisfies the
orthogonality relations \reff{def.orthogonality}
with respect to any family of (not necessarily positive)
measures $\bmu' = (\mu'_1,\ldots,\mu'_r)$
where $\mu'_i$ is any linear combination of $\{ \mu_j \colon n_j \ge n_i \}$.
In particular, if $n_1 \ge n_2 \ge \ldots \ge n_r$,
then we can take $\mu'_i = \sum_{j=1}^i c_{ij} \mu_j$
for any lower-triangular matrix $C = (c_{ij})_{1 \le i,j \le r}$.
That is, $\mu'_i$ is an arbitrary linear combination of
$\mu_1,\ldots,\mu_i$.
This observation will play an important role in what follows
(see Section~\ref{sec.orthogonal} ff.).

The collection $(P_\bfn(x))_{\bfn \in \N^r}$ of
(monic) multiple orthogonal polynomials
of type~II satisfies a variety of recurrences, generalizing the well-known
three-term recurrence for conventional orthogonal polynomials.
Here is one \cite[Theorem~23.1.7 {\em et seq.}\/]{Ismail_05}:
We denote by $\bfe_k$ the multi-index with entry 1 in position $k$
and 0 elsewhere.
For a permutation $\pi$ of $\{1,\ldots,r\}$, we write
$\bfs^{(\pi)}_j = \bfe_{\pi(1)} + \bfe_{\pi(2)} + \ldots + \bfe_{\pi(j)}$
for $1 \le j \le r$.
Then there exist real numbers $a_{\bfn,0}^{(k)}$
($\bfn \in \N^r$, $1 \le k \le r$)
and $a_{\bfn,j}^{(\pi)}$
($\bfn \in \N^r$, $\pi \in \Sym_r$, $1 \le j \le r$) such that
\be
   x \, P_\bfn(x)
   \;=\;
   P_{\bfn + \bfe_k}(x)
   \:+\:
   a_{\bfn,0}^{(k)} \, P_\bfn(x)
   \:+\:
   \sum_{j=1}^r a_{\bfn,j}^{(\pi)} \, P_{\bfn - \bfs^{(\pi)}_j}(x)
   \;,
 \label{eq.MOP.gen_recurrence}
\ee
with the convention that $P_{\bfn - \bfs^{(\pi)}_j}(x) = 0$
whenever one or more of the entries in $\bfn - \bfs^{(\pi)}_j$ is negative.
(Note that the coefficients $a_{\bfn,j}^{(\pi)}$ do not depend on $k$,
 and the coefficients $a_{\bfn,0}^{(k)}$ do not depend on $\pi$.
 But we will never use this fact.)

Now let $j_1,j_2,\ldots$ be an infinite sequence
of elements of $\{1,\ldots,r\}$,
and define a sequence $(\bfn_k)_{k \ge 0}$ of multi-indices in $\N^r$
by $\bfn_k = \sum_{i=1}^k \bfe_{j_i}$.
These multi-indices satisfy $|\bfn_k| = k$
and describe an increasing nearest-neighbor path in $\N^r$
in which the $i$th step is along direction $j_i$.
Now let $\Phat_k(x) \eqdef P_{\bfn_k}(x)$
be the multiple orthogonal polynomial of type~II along this path in $\N^r$.
It then follows from \reff{eq.MOP.gen_recurrence}
that the singly-indexed sequence $(\Phat_n(x))_{n \ge 0}$
satisfies an $(r+2)$-term recurrence of the form
\be
   x \Phat_n(x)   \;=\;  \sum_{k=n-r}^{n+1} \pi_{nk} \, \Phat_k(x)
 \label{eq.prop.Ptildenx.recurrence.2}
\ee
where $\pi_{n,n+1} = 1$ and $\pi_{nk} = 0$ for $k < 0$;
of course the coefficients $\pi_{nk}$ depend
on the choice of nearest-neighbor path.
The recurrence \reff{eq.prop.Ptildenx.recurrence.2}
will play a central role in the remainder of this paper.
Please observe that the coefficients $\pi_{nk}$ in this recurrence
can be collected into a matrix $\Pi = (\pi_{nk})_{n,k \ge 0}$
that is unit-lower-Hessenberg and $(r,1)$-banded:
that is, $\pi_{nk} = 0$ if $k > n+1$ or $k < n-r$, and $\pi_{n,n+1} = 1$.

A particularly important role is played by
the multi-indices $\bfn = (n_1,\ldots,n_r)$ lying on the \textbfit{stepline}:
this is the near-diagonal sequence starting at $(0,0,\ldots,0)$
and following the path
$
   (n,n,\ldots,n) \to
   (n+1,n,\ldots,n) \to
   (n+1,n+1,\ldots,n) \to \ldots \to
   (n+1,n+1,\ldots,n+1) \to \ldots
   \;.
$
In~other words, we define a singly-indexed sequence $(\Ptilde_n(x))_{n \ge 0}$
by
\be
   \Ptilde_n(x)  \;=\;  P_{(n_1,\ldots,n_r)}(x)
   \quad\hbox{where }
   n_i \:=\: \Bigl\lfloor {n+r-i \over r} \Bigr\rfloor
   \;\;\hbox{for } 1 \le i \le r
   \;.
 \label{def.Ptilde.stepline}
\ee
The stepline polynomials $(\Ptilde_n(x))_{n \ge 0}$
are a special case of the nearest-neighbor-path polynomials
$(\Phat_n(x))_{n \ge 0}$,
so they satisfy an $(r+2)$-term recurrence
of the form \reff{eq.prop.Ptildenx.recurrence.2}.

%

\subsection{Production matrices}  \label{subsec.prodmat_intro}

In this subsection we give a brief introduction to the theory of
production matrices \cite{Deutsch_05,Deutsch_09};
see also \cite{Sokal_totalpos} \cite[sections~8.1 and 9.2]{latpath_SRTR}
for further discussion.
In the general theory, the production matrix can be any row-finite
or column-finite matrix.  Here, however, we shall restrict attention
to production matrices that are unit-lower-Hessenberg.
Also, in the general theory the production matrix is usually called $P$;
but here we shall call it $\Pi$ in order to avoid confusion with the
sequence of polynomials $P_n(x)$.

So let $\Pi = (\pi_{ij})_{i,j \ge 0}$ be a unit-lower-Hessenberg matrix
(indexed by $\N$) with entries in a commutative ring $R$:
that is, $\pi_{i,i+1} = 1$ and $\pi_{ij} = 0$ for $j > i+1$.
Then let $A = (a_{nk})_{n,k \ge 0}$ be the matrix defined by
$a_{nk} = (\Pi^n)_{0k}$.
It is easy to see that $A$ is unit-lower-triangular,
i.e.\ $a_{nn} = 1$ and $a_{nk} = 0$ for $k > n$.
Writing out the matrix multiplications explicitly, we have
\be
   a_{nk}
   \;=\; 
   \sum_{i_1,\ldots,i_{n-1}}
      \pi_{0 i_1} \, \pi_{i_1 i_2} \, \pi_{i_2 i_3} \,\cdots\,
        \pi_{i_{n-2} i_{n-1}} \, \pi_{i_{n-1} k}
   \;,
 \label{def.iteration.walk}
\ee
so that $a_{nk}$ is the total weight for all $n$-step walks in $\N$
from $i_0 = 0$ to $i_n = k$, in~which the weight of a walk is the
product of the weights of its steps, and a step from $i$ to $j$
gets a weight $\pi_{ij}$.
(Since $\Pi$ is lower-Hessenberg, these are \textbfit{\L{}ukasiewicz walks},
 i.e.\ the allowed steps are $i \to j$ with $0 \le j \le i+1$.)
Yet another equivalent formulation is to define the entries $a_{nk}$
by the recurrence
\be
   a_{nk}  \;=\;  \sum_{i=0}^\infty a_{n-1,i} \, \pi_{ik}
   \qquad\hbox{for $n \ge 1$}
 \label{def.iteration.bis}
\ee
with the initial condition $a_{0k} = \delta_{0k}$.
We shall call $\Pi$ the \textbfit{production matrix}
and $A$ the \textbfit{output matrix}, and we write $A = \scro(\Pi)$.

These definitions can be given a compact matrix formulation.
Let $\Delta$ be the matrix with 1 on the superdiagonal and 0 elsewhere,
i.e.\ $\Delta_{n,n+1} = 1$ and $\Delta_{nk} = 0$ for $k \neq n+1$
(of course it is unit-lower-Hessenberg).
Then for any matrix $M$ with rows indexed by $\N$,
the product $\Delta M$ is simply $M$ with its zeroth row removed
and all other rows shifted upwards.
(Some authors use the notation $\overline{M} \eqdef \Delta M$.)
The recurrence \reff{def.iteration.bis} can then be written as
\be
   \Delta A  \;=\;  A \Pi
   \;.
 \label{def.iteration.bis.matrixform}
\ee
Since $A$ is unit-lower-triangular, it is invertible,
so \reff{def.iteration.bis.matrixform} is equivalent to
\be
   \Pi  \;=\;  A^{-1} \Delta A
   \;.
 \label{def.iteration.bis.matrixform.2}
\ee
Conversely, it is not difficult to see that
for any unit-lower-triangular matrix $A$,
the matrix $A^{-1} \Delta A$ is unit-lower-Hessenberg.
It therefore follows that for each unit-lower-triangular matrix $A$,
there is a unique unit-lower-Hessenberg matrix $\Pi$
such that $A = \scro(\Pi)$,
and it is given by $\Pi = A^{-1} \Delta A$.

Now let $B = A^{-1}$ be the inverse of $A$
(which is of course also unit-lower-triangular).
We then have a one-to-one correspondence between
unit-lower-Hessenberg matrices~$\Pi$,
unit-lower-triangular matrices~$A$
and unit-lower-triangular matrices~$B$, defined by
\be
   A \:=\: \scro(\Pi) \:=\: B^{-1} \:,\quad
   B \:=\: \scro(\Pi)^{-1} \:=\: A^{-1} \:,\quad
   \Pi \:=\: A^{-1} \Delta A \:=\: B \Delta B^{-1}
   \;.
 \label{eq.A.B.Pi}
\ee
In Section~\ref{sec.prodmat} we will see how the matrices $A$, $B$ and $\Pi$
arise in different characterizations of sequences of monic polynomials.

\bigskip

{\bf Remarks.}
1.  Production matrices are nowadays widely used in enumerative combinatorics:
thus, for instance, the entry for a triangular array
in the On-Line Encyclopedia of Integer Sequences \cite{OEIS}
often gives its production matrix.

2.  Several subclasses of lower-Hessenberg production matrices
are of especial combinatorial interest:
\begin{itemize}
   \item Tridiagonal production matrices correspond to \textbfit{Motzkin walks}
      (i.e.\ the allowed steps are $i \to j$ with $j \in \{i-1,i,i+1\}$)
      and thence to classical J-fractions \cite{Flajolet_80},
      as will be explained in Section~\ref{subsec.classical.CF}.
   \item Toeplitz lower-Hessenberg production matrices generate
      Bell-subgroup \textbfit{Riordan arrays};  and more generally,
      1-almost-Toeplitz lower-Hessenberg production matrices
      (i.e.\ lower-Hessenberg matrices that are Toeplitz
       except for the zeroth column) generate Riordan arrays.
      See \cite{Shapiro_91,Sprugnoli_94,Barry_16,Shapiro_22}
      for introductions to Riordan arrays,
      and \cite{Deutsch_05,He_15,Sokal_totalpos} for the just-quoted theorems
      on their production matrices.
   \item Lower-Hessenberg production matrices of the form
      $\pi_{nk} = (n!/k!) \, (z_{n-k} + k a_{n-k+1})$
      generate \textbfit{exponential Riordan arrays}
      \cite[pp.~217--218]{Barry_16} \cite[Theorem~8.2]{latpath_lah}.
\end{itemize}

3.  When the commutative ring $R$ is equipped with a partial order,
production matrices are also a powerful tool for attacking problems
related to total positivity \cite{Sokal_totalpos}.
In particular, the total positivity of the production matrix $\Pi$
is a sufficient (but far~from necessary) condition for the
total positivity of its output matrix $\scro(\Pi)$
and for the Hankel-total positivity of the zeroth-column sequence
of $\scro(\Pi)$.
See \cite{Sokal_totalpos} \cite[sections~8.1 and 9.2]{latpath_SRTR}
    \cite{latpath_lah,latpath_laguerre,forests_totalpos}
for precise statements, proofs, and further discussion
and applications.
\myendremark

\subsection{Classical continued fractions}   \label{subsec.classical.CF}

As preparation for the discussion of branched continued fractions
in Section~\ref{subsec.BCF},
as well as for some applications later in this paper,
we begin by giving a very brief review of selected aspects of the theory
of classical continued fractions (J-fractions and S-fractions).
We will follow the notation and terminology used nowadays by
combinatorialists \cite{Flajolet_80},
as this is the most appropriate for our work;
but we will also point out the translation to the formalism
employed in the classical analysis literature
on continued fractions
\cite{Perron,Wall_48,Jones_80,Lorentzen_92,Cuyt_08}
and the moment problem
\cite{Stieltjes_1894,Shohat_43,Ahiezer_62,Akhiezer_65,Schmudgen_17}.

We shall consider continued fractions of either Stieltjes (S) type,
\be
   f(t)
   \;=\;
   \sum_{n=0}^\infty a_n t^n
   \;=\;
   \cfrac{1}{1 - \cfrac{\alpha_1 t}{1 - \cfrac{\alpha_2 t}{1 - \cdots}}}
   \label{def.Stype}
   \;\;,
\ee
or Jacobi (J) type,
\be
   f(t)
   \;=\;
   \sum_{n=0}^\infty a_n t^n
   \;=\;
   \cfrac{1}{1 - \gamma_0 t - \cfrac{\beta_1 t^2}{1 - \gamma_1 t - \cfrac{\beta_2 t^2}{1 - \cdots}}}
   \label{def.Jtype}
   \;\;.
\ee
Here these expressions are to be interpreted as
formal power series in the indeterminate $t$;
we do {\em not}\/ wish to address questions of convergence.
Thus, the continued-fraction coefficients $\balpha = (\alpha_n)_{n \ge 1}$,
$\bbeta = (\beta_n)_{n \ge 1}$ and $\bgamma = (\gamma_n)_{n \ge 0}$
are sequences in a commutative ring $R$,
and they determine the sequence $\ba = (a_n)_{n \ge 0}$ of Taylor coefficients
by formal expansion of the continued fraction.
Indeed, it is conceptually simplest to consider
$\balpha,\bbeta,\bgamma$ as algebraic indeterminates;
then the $a_n$ are polynomials with integer coefficients
in these indeterminates:
\begin{eqnarray}
   \sum_{n=0}^\infty S_n(\balpha) \, t^n
   & = &
   \cfrac{1}{1 - \cfrac{\alpha_1 t}{1 - \cfrac{\alpha_2 t}{1 - \cdots}}}
   \label{def.SRpoly}
       \\[4mm]
   \sum_{n=0}^\infty J_n(\bbeta,\bgamma) \, t^n
   & = &
   \cfrac{1}{1 - \gamma_0 t - \cfrac{\beta_1 t^2}{1 - \gamma_1 t - \cfrac{\beta_2 t^2}{1 - \cdots}}}
   \label{def.JRpoly}
\end{eqnarray}
We call $S_n(\balpha)$ the \textbfit{Stieltjes--Rogers polynomials},
and $J_n(\bbeta,\bgamma)$ the \textbfit{Jacobi--Rogers polynomials}.

In a seminal 1980 paper, Flajolet \cite{Flajolet_80}
gave a combinatorial interpretation of the Stieltjes--Rogers
and Jacobi--Rogers polynomials in terms of lattice paths.
We recall that a \textbfit{Motzkin path} of length $n$
is a path in the right quadrant $\N \times \N$,
starting at $(0,0)$ and ending at $(n,0)$,
using steps $(1,1)$ [``rise''], $(1,0)$ [``level step'']
and $(1,-1)$ [``fall''].
More generally, a \textbfit{Motzkin path at level $\bm{k}$}
is a path in $\N \times \N_{\ge k}$,
starting at $(0,k)$ and ending at $(n,k)$,
using the same steps.
A Motzkin path is called a \textbfit{Dyck path} if it has no level steps;
obviously a Dyck path must have even length.

\begin{theorem}[Flajolet \protect\cite{Flajolet_80}]
   \label{thm.flajolet}
\quad\hfill\break
\vspace*{-6mm}
\begin{itemize}
   \item[(a)]  The Jacobi--Rogers polynomial $J_n(\bbeta,\bgamma)$
is the generating polynomial for Motzkin paths of length $n$,
in which each rise gets weight~1,
each level step at height $i$ gets weight $\gamma_i$,
and each fall from height $i$ gets weight $\beta_i$.
   \item[(b)]  The Stieltjes--Rogers polynomial $S_n(\balpha)$
is the generating polynomial for Dyck paths of length $2n$,
in which each rise gets weight~1
and each fall from height $i$ gets weight $\alpha_i$.
\end{itemize}
\end{theorem}

\par\medskip\noindent{\sc Proof} \cite{Flajolet_80}.
(a) For each $k \ge 0$, let $f_k(t)$ be the generating function
for Motzkin paths at level $k$ (of arbitrary length)
in which each rise gets weight~1,
each level step at height $i$ gets weight $\gamma_i$,
each fall from height $i$ gets weight $\beta_i$,
and each step of any kind gets an additional weight $t$.
It is a formal power series in the indeterminate $t$,
with coefficients that are polynomials in $\bbeta,\bgamma$.

Now let $\scrp$ be any Motzkin path at level $k$;
and if it is of nonzero length,
split it at its first return to height $k$,
yielding $\scrp = \scrp' \, \scrp''$.
Then $\scrp'$ is either a single level step at height $k$,
or else a path of the form $U \scrp_{k+1} D$
where $U$ is a rise $k \to k+1$,
$\scrp_{k+1}$ is an arbitrary Motzkin path at level $k+1$,
and $D$ is a fall $k+1 \to k$.
Furthermore, $\scrp''$ is an arbitrary Motzkin path at level $k$.
We thus deduce the functional equation
\be
   f_k(t)
   \;=\;
   1 \:+\: \gamma_k t \,f_k(t) \:+\: \beta_{k+1} t^2 \, f_{k+1}(t) \, f_k(t)
 \label{eq.JRfk.1}
\ee
or equivalently
\be
   f_k(t)
   \;=\;
   {1 \over  1 \:-\: \gamma_k t \:-\: \beta_{k+1} t^2 \, f_{k+1}(t)}
   \;.
 \label{eq.JRfk.2}
\ee
Iterating \reff{eq.JRfk.2}, we see immediately that $f_k$
is given by the continued fraction
\be
   f_k(t)
   \;=\;
   \cfrac{1}{1 - \gamma_k t - \cfrac{\beta_{k+1} t^2}{1 - \gamma_{k+1} t - \cfrac{\beta_{k+2} t^2}{1 - \cdots}}}
\ee
and in particular that $f_0$ is given by \reff{def.JRpoly}.

(b) This follows from part~(a) by setting $\bgamma = \bzero$,
renaming $\bbeta$ as $\balpha$, and renaming $t^2$ as $t$.
\qed

{\bf Remarks.}
1.  In the function-theoretic literature on the moment problem
\cite{Stieltjes_1894,Shohat_43,Ahiezer_62,Akhiezer_65,Schmudgen_17}
and continued fractions
\cite{Perron,Wall_48,Jones_80,Lorentzen_92,Cuyt_08},
the generating function for a sequence $\ba = (a_n)_{n \ge 0}$
of real numbers is most often written in the form
\be
   F(z)  \;=\;  {1 \over z} f\Bigl( {1 \over z} \Bigr)
         \;=\; \sum_{n=0}^\infty {a_n \over z^{n+1}}
   \;.
 \label{def.Fz}
\ee
This formulation has the property that if $\ba$ is a moment sequence
with representing measure $\mu$,
i.e.\ $a_n = \int_{-\infty}^\infty x^n \, d\mu(x)$,
then the Stieltjes transform
\be
   F(z)  \;\eqdef\; \int\limits_{-\infty}^\infty {d\mu(x) \over z-x}
\ee
is analytic in the upper half-plane $\Im z > 0$
and has the series \reff{def.Fz} as its large-$z$ asymptotic expansion,
uniformly in each sector $\epsilon \le |\arg z| \le \pi-\epsilon$
\cite[p.~27]{Shohat_43}.

Given a power series of the form \reff{def.Fz},
the S-type continued fraction is then written in the form
\cite[p.~viii]{Shohat_43} \cite[p.~329]{Wall_48}
\be
   F(z)  \;=\;
   \cfrac{1}{l_1 z - \cfrac{1}{l_2 - \cfrac{1}{l_3 z -  \cfrac{1}{l_4 - \cdots}}}}
   \label{def.Stype.F}
\ee
(note that $l_n$ is multiplied by $z$ for $n$ odd but not for $n$ even),
which is easily seen to be equivalent to \reff{def.Stype}
if we normalize to $a_0 = 1$ (hence $l_1 = 1$)
and then set $\alpha_1 = 1/l_2$ and
$\alpha_n = 1/(l_n l_{n+1})$ for $n \ge 2$;
the reverse translation is
\begin{subeqnarray}
   l_{2k-1}  & = & {\alpha_1 \alpha_3 \cdots \alpha_{2k-3}
                    \over
                    \alpha_2 \alpha_4 \cdots \alpha_{2k-2}
                   }
       \\[2mm]
   l_{2k}    & = & {\alpha_2 \alpha_4 \cdots \alpha_{2k-2}
                    \over
                    \alpha_1 \alpha_3 \cdots \alpha_{2k-1}
                   }
 \label{eq.alpha-to-l}
\end{subeqnarray}
%
Likewise, the J-type continued fraction is written in the form
\cite[pp.~viii, 31]{Shohat_43}
\be
   F(z)  \;=\;
   \cfrac{\lambda_1}{z - c_1 - \cfrac{\lambda_2}{z - c_2 - \cfrac{\lambda_3}{z - c_3 - \cfrac{\lambda_4}{z - c_4 - \cdots}}}}
   \label{def.Jtype.F}
   \;\;,
\ee
which is easily seen to be equivalent to \reff{def.Jtype}
if we normalize to $\lambda_1 = 1$
and then set $\gamma_n = c_{n+1}$ and $\beta_n = \lambda_{n+1}$.

2.  My use of the terms ``S-fraction'' and ``J-fraction'' follows
the general practice in the combinatorial literature,
starting with Flajolet \cite{Flajolet_80}.
The classical literature on continued fractions
\cite{Perron,Wall_48,Jones_80,Lorentzen_92,Cuyt_08}
generally uses a different terminology.
For instance, Jones and Thron \cite[pp.~128--129, 386--389]{Jones_80}
use the term ``regular C-fraction''
for (a minor variant of) what I~have called an S-fraction;
they call it an ``S-fraction'' if all $\alpha_n < 0$.
They use the term ``associated continued fraction''
for (a minor variant of) what I~have called a J-fraction,
and use the term ``J-fraction'' for \reff{def.Jtype.F}
with $\lambda_n \neq 0$.

3.  It is worth observing that an S-fraction can always be transformed
into a J-fraction by {\em contraction}\/
\cite[p.~21]{Wall_48} \cite[p.~V-31]{Viennot_83}:
namely, \reff{def.Stype} and \reff{def.Jtype} are equal if
\begin{subeqnarray}
   \gamma_0  & = &  \alpha_1
       \slabel{eq.contraction_even.coeffs.a}   \\
   \gamma_n  & = &  \alpha_{2n} + \alpha_{2n+1}  \qquad\hbox{for $n \ge 1$}
       \slabel{eq.contraction_even.coeffs.b}   \\
   \beta_n  & = &  \alpha_{2n-1} \alpha_{2n}
       \slabel{eq.contraction_even.coeffs.c}
 \label{eq.contraction_even.coeffs}
\end{subeqnarray}
See \cite[pp.~20--22]{Wall_48} for the classic algebraic proof;
see \cite[Lemmas~1 and 2]{Dumont_94b}
\cite[proof of Lemma~1]{Dumont_95} \cite[Lemma~4.5]{DiFrancesco_10}
for a very simple variant algebraic proof;
and see \cite[pp.~{V-31}--{V-32}]{Viennot_83}
for an enlightening combinatorial proof,
based on defining a Motzkin path by grouping pairs of steps in a Dyck path.
The reverse transformation --- from J-fraction to S-fraction ---
is generically possible if the coefficient ring $R$ is a field,
but not in general otherwise.
\myendremark

\bigskip

Let us now generalize these definitions;
we concentrate on the case of J-fractions,
but similar constructions can be applied to S-fractions.
A \textbfit{partial Motzkin path} of length $n$
is a path in the right quadrant $\N \times \N$,
starting at $(0,0)$ and ending at some point $(n,k)$,
using the same steps as before.
Let $J_{n,k}(\bbeta,\bgamma)$ be the generating polynomial
for partial Motzkin paths from $(0,0)$ to $(n,k)$,
in which each rise gets weight~1,
each level step at height $i$ gets weight $\gamma_i$,
and each fall from height $i$ gets weight $\beta_i$.
We therefore have an infinite unit-lower-triangular array
$\sfJ = \big( J_{n,k}(\bbeta,\bgamma) \big)_{n,k \ge 0}$
in which the first ($k=0$) column displays
the ordinary Jacobi--Rogers polynomials $J_{n,0} = J_n$.
It is immediate from the definition of $J_{n,k}$
that the matrix $\sfJ$ is the output matrix $\scro(\Pi)$
corresponding to the tridiagonal production matrix
\be
   \Pi
   \;=\;
   \begin{bmatrix}
      \gamma_0   & 1           &            &        &         \\
      \beta_1    & \gamma_1    & 1          &        &         \\
                 & \beta_2     & \gamma_2   & 1      &         \\
                 &             & \ddots     & \ddots & \ddots
   \end{bmatrix}
 \label{def.Jnl.production.bis.0}
\ee
that generates Motzkin walks with the given weights
[cf.\ \reff{def.iteration.walk}].
We then have the following beautiful fact:

\begin{proposition}[$LDL^{\rm T}$ factorization of the Hankel matrix of Jacobi--Rogers polynomials]
   \label{prop.LDLT}
The Hankel matrix of Jacobi--Rogers polynomials,
\be
   H_\infty(\bJ)  \;\eqdef\; \big( J_{n+n'}(\bbeta,\bgamma) \big)_{n,n' \ge 0}
   \;,
\ee
has the factorization
\be
   H_\infty(\bJ)  \;=\;  \sfJ D \sfJ^{\rm T}
  \label{eq.jacobi.LDLT.0}
\ee
where $D = \diag(1,\beta_1, \beta_1 \beta_2, \ldots)$
is the diagonal matrix with entries
\be
   D_{kk}  \;=\;  \prod_{i=1}^{k} \beta_i
\ee
for $k \ge 0$.
\end{proposition}

\proof
It suffices to note the identity
\cite[p.~351]{Aigner_01a} \cite[Remark~2.2]{Ismail_10}
\be
   J_{n+n',0}(\bbeta,\bgamma)
   \;=\;
   \sum_{\ell=0}^\infty
   J_{n,\ell}(\bbeta,\bgamma)
   \biggl( \prod_{i=1}^\ell \beta_i \! \biggr)
   J_{n',\ell}(\bbeta,\bgamma)
   \;,
 \label{eq.Jnl.hankel.0}
\ee
which arises from splitting a Motzkin path of length $n+n'$
into its first $n$ steps and its last $n'$ steps,
and then imagining the second part run backwards:
the factor $\prod_{i=1}^k \beta_i$ arises from the fact that
when we reversed the path we interchanged rises with falls
and thus lost a factor $\prod_{i=1}^k \beta_i$
for those falls that were not paired with rises.
The identity \reff{eq.Jnl.hankel.0} can be written in matrix form as
\reff{eq.jacobi.LDLT.0}.
\qed

{\bf Remarks.}
1. The reversal argument employed in this proof
can be rewritten purely algebraically as follows:
Note first that the tridiagonal matrix \reff{def.Jnl.production.bis.0}
satisfies
\be
   \Pi^{\rm T}  \;=\;  D^{-1} \Pi D
 \label{eq.Pi.diagsim}
\ee
where $D = \diag(1,\beta_1, \beta_1 \beta_2, \ldots)$.
(Here we work in the ring $\Z[\bbeta,\bbeta^{-1},\bgamma]$
 of Laurent polynomials in $\bbeta$.)
On the other hand, it is a general fact \cite{Sokal_totalpos}
that if $\Pi$ is a production matrix
and $\scroo_0(\Pi)$ is the zeroth-column sequence of $\scro(\Pi)$,
then
\be
   H_\infty(\scroo_0(\Pi))
   \;=\;
   \scro(\Pi) \, {\scro(\Pi^{\rm T})}^{\rm T}
   \;.
\ee
And finally, it is a general fact \cite{Sokal_totalpos}
that if $M$ is an invertible lower-triangular matrix
satisfying $M_{00} = 1$, then $\scro(M^{-1} \Pi M) = \scro(\Pi) \, M$.
Putting all this together, we have
\begin{subeqnarray}
   H_\infty(\scroo_0(\Pi))
   & = &
   \scro(\Pi) \, {\scro(\Pi^{\rm T})}^{\rm T}
         \\[2mm]
   & = &
   \scro(\Pi) \, {\scro(D^{-1} \Pi D)}^{\rm T}
         \\[2mm]
   & = &
   \scro(\Pi) \, [\scro(\Pi) \, D]^{\rm T}
         \\[2mm]
   & = &
   \scro(\Pi) \, D \, \scro(\Pi)^{\rm T}
      \;,
\end{subeqnarray}
as asserted in Proposition~\ref{prop.LDLT}.
Obviously, this proof relies crucially on the fact that
the reversal of a Motzkin path is again a Motzkin path,
or equivalently on the fact that
the production matrix $\Pi$ is tridiagonal
and therefore symmetric up to a diagonal similarity transformation
\reff{eq.Pi.diagsim}.

2. The factorization \reff{eq.jacobi.LDLT.0}
was found more than a century ago by
Stieltjes \cite{Stieltjes_1889,Stieltjes_1894},
albeit without the interpretation in terms of Motzkin paths.
More precisely, Stieltjes
\cite{Stieltjes_1889} \cite[pp.~J.18--J.19]{Stieltjes_1894}
found the analogous factorization for S-fractions.
The factorization for J-fractions can be found in Wall's 1948 book
\cite[Theorem~53.1]{Wall_48}, among other places.
\myendremark

Taking the determinant of the $n \times n$ leading principal submatrix
on both sides of \reff{eq.jacobi.LDLT.0}, we obtain a classical formula
\cite[Theorem~51.1]{Wall_48} 
for the Hankel determinants of a J-fraction:

\begin{corollary}[Hankel determinants of the Jacobi--Rogers polynomials]
Let $\Delta_n = \det \big( J_{i+j}(\bbeta,\bgamma) \big)_{0 \le i,j \le n-1}$
be the $n \times n$ leading principal minor of the Hankel matrix
$H_\infty(\bJ)$.  Then
\be
   \Delta_{n+1}  \;=\;  \beta_1^n \beta_2^{n-1} \cdots \beta_{n-1}^2 \beta_n
   \;.
 \label{eq.deltas.Jfrac.0}
\ee
\end{corollary}

In particular, if $R = \R$ and all $\beta_i > 0$,
then the Hankel matrix $H_\infty(\bJ)$ is positive-definite,
which implies that the underlying sequence $(J_n(\bbeta,\bgamma))_{n \ge 0}$
is a Hamburger moment sequence
with a representing measure of infinite support
\cite{Shohat_43,Ahiezer_62,Akhiezer_65,Schmudgen_17}.

\subsection{Branched continued fractions}   \label{subsec.BCF}

In this subsection we give a very brief introduction to the theory of
branched continued fractions,
limiting attention for simplicity to branched S-fractions;
our treatment follows \cite{latpath_SRTR},
where many more details and applications can be found.

Fix an integer $m \ge 1$.
An \textbfit{$\bm{m}$-Dyck path}
\cite{Aval_08,Cameron_16,Prodinger_16,latpath_SRTR}
is a path in the upper half-plane $\Z \times \N$,
starting and ending on the horizontal axis,
using steps $(1,1)$ [``rise'' or ``up step'']
and $(1,-m)$ [``$m$-fall'' or ``down step''].
More generally, an \textbfit{$\bm{m}$-Dyck path at level $\bm{k}$}
is a path in $\Z \times \N_{\ge k}$,
starting and ending at height $k$,
using steps $(1,1)$ and $(1,-m)$.
Since the number of up steps must equal $m$ times the number of down steps,
the length of an $m$-Dyck path must be a multiple of $m+1$.

Now let $\balpha = (\alpha_i)_{i \ge m}$ be an infinite set of indeterminates.
Then \cite{latpath_SRTR}
the \textbfit{$\bm{m}$-Stieltjes--Rogers polynomial} of order~$n$,
denoted $S^{(m)}_n(\balpha)$, is the generating polynomial
for $m$-Dyck paths of length~$(m+1)n$ in which each rise gets weight~1
and each $m$-fall from height~$i$ gets weight $\alpha_i$.
Clearly $S_n^{(m)}(\balpha)$ is a homogeneous polynomial
of degree~$n$ with nonnegative integer coefficients.

Let $f_0(t) = \sum_{n=0}^\infty S^{(m)}_n(\balpha) \, t^n$
be the ordinary generating function for $m$-Dyck paths with these weights;
and more generally, let $f_k(t)$ be the ordinary generating function
for $m$-Dyck paths at level $k$ with these same weights.
(Obviously $f_k$ is just $f_0$ with each $\alpha_i$ replaced by $\alpha_{i+k}$;
 but we shall not explicitly use this fact.)
Then straightforward combinatorial arguments \cite[Section~2.3]{latpath_SRTR},
similar to those used in the proof of Theorem~\ref{thm.flajolet},
lead to the functional equation
\be
   f_k(t)  \;=\;  1 \:+\: \alpha_{k+m} t \, f_k(t) \, f_{k+1}(t) \,\cdots\, f_{k+m}(t)
 \label{eq.mSRfk.1}
\ee
or equivalently
\be
   f_k(t)  \;=\;  {1 \over 1 \:-\: \alpha_{k+m} t \, f_{k+1}(t) \,\cdots\, f_{k+m}(t)}
   \;.
 \label{eq.mSRfk.2}
\ee
Iterating \reff{eq.mSRfk.2}, we see immediately that $f_k$
is given by the {\em branched continued fraction}\/
\begin{subeqnarray}
   f_k(t)
   & = &
   \cfrac{1}
         {1 \,-\, \alpha_{k+m} t
            \prod\limits_{i_1=1}^{m}
                 \cfrac{1}
            {1 \,-\, \alpha_{k+m+i_1} t
               \prod\limits_{i_2=1}^{m}
               \cfrac{1}
            {1 \,-\, \alpha_{k+m+i_1+i_2} t
               \prod\limits_{i_3=1}^{m}
               \cfrac{1}{1 - \cdots}
            }
           }
         }
%
      \slabel{eq.fk.mSfrac.a} \\[2mm]
   & = &
\Scale[0.6]{
   \cfrac{1}{1 - \cfrac{\alpha_{k+m} t}{
     \Biggl( 1 - \cfrac{\alpha_{k+m+1} t}{
        \Bigl( 1  - \cfrac{\alpha_{k+m+2} t}{(\cdots) \,\cdots\, (\cdots)} \Bigr)
        \,\cdots\,
        \Bigl( 1  - \cfrac{\alpha_{k+2m+1} t}{(\cdots) \,\cdots\, (\cdots)} \Bigr)
       }
     \Biggr)
     \,\cdots\,
     \Biggl( 1 - \cfrac{\alpha_{k+2m} t}{
        \Bigl( 1  - \cfrac{\alpha_{k+2m+1} t}{(\cdots) \,\cdots\, (\cdots)} \Bigr)
        \,\cdots\,
        \Bigl( 1  - \cfrac{\alpha_{k+3m} t}{(\cdots) \,\cdots\, (\cdots)} \Bigr)
       }
     \Biggr)
    }
   }
}
     \nonumber \\
 \slabel{eq.fk.mSfrac.b}
 \label{eq.fk.mSfrac}
\end{subeqnarray}
and in particular that $f_0$ is given by
the specialization of \reff{eq.fk.mSfrac} to $k=0$.
We shall call the right-hand side of \reff{eq.fk.mSfrac}
an \textbfit{$\bm{m}$-branched Stieltjes-type continued fraction},
or \textbfit{$\bm{m}$-branched S-fraction} for short.

\medskip

{\bf Remark.}
In truth, we hardly ever use the branched continued fraction
\reff{eq.fk.mSfrac};
instead, we work directly with the $m$-Dyck paths
and/or with the recurrence \reff{eq.mSRfk.1}/\reff{eq.mSRfk.2}
that their generating functions satisfy.
\myendremark

\medskip

We now generalize these definitions as follows.
A \textbfit{partial $\bm{m}$-Dyck path}
is a path in the upper half-plane $\Z \times \N$,
starting on the horizontal axis but ending anywhere,
using steps $(1,1)$ [``rise'']
and $(1,-m)$ [``$m$-fall''].
A partial $m$-Dyck path starting at $(0,0)$
must stay always within the set
$V_m = \{ (x,y) \in \Z \times \N \colon\: x=y \bmod m+1 \}$.

Now let $\balpha = (\alpha_i)_{i \ge m}$ be an infinite set of indeterminates,
and let $S^{(m)}_{n,k}(\balpha)$ be the generating polynomial
for partial $m$-Dyck paths from $(0,0)$ to ${((m+1)n,(m+1)k)}$
in~which each rise gets weight~1
and each $m$-fall from height~$i$ gets weight $\alpha_i$.
We call the $S^{(m)}_{n,k}$ the
\textbfit{generalized $\bm{m}$-Stieltjes--Rogers polynomials}.
Obviously $S^{(m)}_{n,k}$ is nonvanishing only for $0 \le k \le n$,
and $S^{(m)}_{n,n} = 1$.
We therefore have an infinite unit-lower-triangular array
$\sfS^{(m)} = \big( S^{(m)}_{n,k}(\balpha) \big)_{n,k \ge 0}$
in which the first ($k=0$) column displays
the ordinary $m$-Stieltjes--Rogers polynomials $S^{(m)}_{n,0} = S^{(m)}_n$.

The production matrix for the triangle $\sfS^{(m)}$
was found in \cite[sections~7.1 and 8.2]{latpath_SRTR}.
We begin by defining some special matrices $M = (m_{ij})_{i,j \ge 0}$:
\begin{itemize}
   \item $L(s_1,s_2,\ldots)$ is the lower-bidiagonal matrix
       with 1 on the diagonal and $s_1,s_2,\ldots$ on the subdiagonal:
\be
   L(s_1,s_2,\ldots)
   \;=\;
   \begin{bmatrix}
      1  &     &     &     &    \\
      s_1 & 1  &     &     &    \\
          & s_2 & 1  &     &    \\
          &     & s_3 & 1  &    \\
          &     &     & \ddots & \ddots
   \end{bmatrix}
   \;.
 \label{def.L}
\ee
   \item $U^\star(s_1,s_2,\ldots)$ is the upper-bidiagonal matrix
       with 1 on the superdiagonal and $s_1,s_2,\ldots$ on the diagonal:
\be
   U^\star(s_1,s_2,\ldots)
   \;=\;
   \begin{bmatrix}
      s_1 & 1   &     &     &     &    \\
          & s_2 & 1   &     &     &    \\
          &     & s_3 & 1   &     &    \\
          &     &     & s_4 & 1   &    \\
          &     &     &     & \ddots & \ddots
   \end{bmatrix}
   \;.
 \label{def.Ustar}
\ee
\end{itemize}
Then the production matrix for the triangle $\sfS^{(m)}$ is
\begin{eqnarray}
   P^{(m)\mathrm{S}}(\balpha)
   & \eqdef &
   L(\alpha_{m+1}, \alpha_{2m+2}, \alpha_{3m+3}, \ldots)
   \:
   L(\alpha_{m+2}, \alpha_{2m+3}, \alpha_{3m+4}, \ldots)
   \:\cdots\:
   \hspace*{1cm}
       \nonumber \\
   & & \qquad
   L(\alpha_{2m}, \alpha_{3m+1}, \alpha_{4m+2}, \ldots)
   \:
   U^\star(\alpha_m, \alpha_{2m+1}, \alpha_{3m+2}, \ldots)
   \;,
   \hspace*{1cm}
 \label{eq.prop.contraction}
\end{eqnarray}
that is, the product of $m$ factors $L$ and one factor $U^\star$
\cite[Proposition~8.2]{latpath_SRTR}.

Let us remark, finally, that there is (as far as I~know) no analogue of
Proposition~\ref{prop.LDLT} for $m$-S-fractions with $m > 1$,
since the reversal of an $m$-Dyck path is not an $m$-Dyck path.

\section{Production matrix for a sequence of monic polynomials}
   \label{sec.prodmat}


In this section we prove some very elementary --- but important ---
relations between sequences of monic polynomials,
the linear recurrences they satisfy,
and their dual sequences of linear functionals.
All of these properties will be re-expressed in a convenient matrix form,
using the theory of production matrices (Section~\ref{subsec.prodmat_intro}).
We conclude this section with some more delicate matters concerning
``expectation values'' of products of polynomials,
culminating in Open Problem~\ref{openproblem.expectation}.

%
%
%
\subsection{Linear functionals}

Let $R[x]$ be the ring of polynomials in one indeterminate~$x$,
with coefficients in $R$;  it is an $R$-module.
(If $R$ is a field, then $R[x]$ is a vector space over $R$.)
A \textbfit{linear functional} (more precisely, an $R$-linear functional)
on $R[x]$ is a map $\scrl \colon R[x] \to R$
satisfying $\scrl(a p(x) + b q(x)) = a \scrl(p(x)) + b \scrl(q(x))$
for all $a,b \in R$ and $p(x), q(x) \in R[x]$.
To each linear functional $\scrl \colon R[x] \to R$
there is naturally associated a sequence $(\ell_n)_{n \ge 0}$
of elements of $R$, defined by $\ell_n = \scrl(x^n)$.
And conversely, to every sequence $(\ell_n)_{n \ge 0}$
there is associated a unique linear functional $\scrl$
satisfying $\scrl(x^n) = \ell_n$,
namely $\scrl(\sum_{n=0}^N c_n x^n) = \sum_{n=0}^N c_n \ell_n$.
We call $(\ell_n)_{n \ge 0}$ the \textbfit{moment sequence}
of the linear functional $\scrl$.

Now let $(\scrl_k)_{k \ge 0}$ be a sequence of such linear functionals.
We form the matrix $A = (a_{nk})_{n,k \ge 0}$ whose columns are the
moment sequences of these linear functionals,
i.e.\ $a_{nk} = \scrl_k(x^n)$.
We call $A$ the \textbfit{moment matrix}
for the sequence $(\scrl_k)_{k \ge 0}$ of linear functionals.
We say that the sequence $(\scrl_k)_{k \ge 0}$ is {\em normalized}\/
in case the matrix $A$ is unit-lower-triangular,
i.e.\ $\scrl_k(x^n) = 0$ for $n < k$ and $\scrl_k(x^k) = 1$.

%
%
\subsection{Sequences of monic polynomials}
By a \textbfit{sequence of monic polynomials}
we mean a sequence $(P_n(x))_{n \ge 0}$ of polynomials
(with coefficients in $R$)
such that $P_n(x)$ has degree $n$ and leading coefficient 1.
We can assemble the coefficients of these polynomials
into a unit-lower-triangular matrix $B = (b_{nk})_{n,k \ge 0}$
by writing $P_n(x) = \sum\limits_{k=0}^n b_{nk} \, x^k$.
There is obviously a one-to-one correspondence between
unit-lower-triangular matrices and sequences of monic polynomials.
We call $B$ the \textbfit{coefficient matrix}
for the sequence $(P_n(x))_{n \ge 0}$ of polynomials.

Now let $A = (a_{nk})_{n,k \ge 0}$ be the inverse matrix to $B$,
i.e.\ $A = B^{-1}$.
Then we obviously have $x^n = \sum\limits_{k=0}^n a_{nk} \, P_k(x)$.

%
%
%
\subsection{Duality}

Let $\bfP = (P_n(x))_{n \ge 0}$ be a sequence of monic polynomials,
and let $\bfL = (\scrl_k)_{k \ge 0}$ be a sequence of linear functionals.
We say that $\bfP$ and $\bfL$ are \textbfit{dual} to each other
in case $\scrl_k(P_n(x)) = \delta_{kn}$ for all $k,n \ge 0$.
The fundamental result concerning such duality is very simple:

\begin{proposition}[Sequence of monic polynomials and its dual sequence
   of linear functionals]
   \label{prop.Pnx.dual}
Given any sequence $(P_n(x))_{n \ge 0}$ of monic polynomials,
there exists a unique sequence $(\scrl_k)_{k \ge 0}$ of
linear functionals that satisfies $\scrl_k(P_n(x)) = \delta_{kn}$,
and it is normalized.

Conversely, given any normalized sequence $(\scrl_k)_{k \ge 0}$ of
linear functionals, there exists a unique sequence $(P_n(x))_{n \ge 0}$
of monic polynomials that satisfies $\scrl_k(P_n(x)) = \delta_{kn}$.

The relation between these sequences is:
The moment matrix $A$ of the sequence $(\scrl_k)_{k \ge 0}$
and the coefficient matrix $B$ of the sequence $(P_n(x))_{n \ge 0}$
are inverses of each other.
\end{proposition}

\proof
Using $P_n(x) = \sum_{j=0}^n b_{nj} \, x^j$ and $\scrl_k(x^j) = a_{jk}$,
we see that the condition $\scrl_k(P_n(x)) = \delta_{kn}$
is equivalent to the matrix equation $BA = I$.
If $(P_n(x))_{n \ge 0}$ is a sequence of monic polynomials,
then $B$ is unit-lower-triangular, and $BA = I$ has the unique solution
$A = B^{-1}$; moreover, $A$ is unit-lower-triangular.
And conversely, if $(\scrl_k)_{k \ge 0}$ is a
normalized sequence of linear functionals,
then $A$ is unit-lower-triangular, and $BA = I$ has the unique solution
$B = A^{-1}$; moreover, $B$ is unit-lower-triangular.
\qed

%
%
%
\subsection{Linear recurrence $\longleftrightarrow$ production matrix}

The next result, which is only slightly more complicated,
connects a sequence $(P_n(x))_{n \ge 0}$ of monic polynomials
with the unique linear recurrence (of a certain standard form) that defines it:

\begin{proposition}[Sequence of monic polynomials and its defining recurrence]
   \label{prop.Pnx.recurrence}
Given any sequence $(P_n(x))_{n \ge 0}$ of monic polynomials,
there exists a unique unit-lower-Hessenberg matrix
$\Pi = (\pi_{nk})_{n,k \ge 0}$ such that
\be
   P_{n+1}(x)
   \;=\;
   (x - \pi_{nn}) \, P_n(x) \:-\: \sum_{k=0}^{n-1} \pi_{nk} \, P_k(x)
 \label{eq.prop.Pnx.recurrence.1}
\ee
or equivalently
\be
   x P_n(x)   \;=\;  \sum_{k=0}^{n+1} \pi_{nk} \, P_k(x)
   \;.
 \label{eq.prop.Pnx.recurrence.2}
\ee

And conversely, given any unit-lower-Hessenberg matrix
$\Pi = (\pi_{nk})_{n,k \ge 0}$,
there exists a unique sequence $(P_n(x))_{n \ge 0}$ of polynomials
satisfying \reff{eq.prop.Pnx.recurrence.1}/\reff{eq.prop.Pnx.recurrence.2}
with the initial condition $P_0(x) = 1$, and it is monic.

The relation between these objects is:
The coefficient matrix $B$ of the sequence $(P_n(x))_{n \ge 0}$
satisfies $B = \scro(\Pi)^{-1}$ or equivalently $\Pi = B \Delta B^{-1}$.
\end{proposition}

\proof
Let $(P_n(x))_{n \ge 0}$ be a sequence of monic polynomials
with coefficient matrix $B$.
Substituting $P_n(x) = \sum_{j=0}^n b_{nj} \, x^j$
into \reff{eq.prop.Pnx.recurrence.2}
and extracting the coefficient of $x^j$,
we see that \reff{eq.prop.Pnx.recurrence.2} is equivalent to
\be
   b_{n,j-1}  \;=\;  (\Pi B)_{nj}
\ee
or in other words
\be
   B \Delta  \;=\;  \Pi B
   \;.
\ee
Since $B$ is unit-lower-triangular and hence invertible,
this equation has the unique solution $\Pi = B \Delta B^{-1}$.

The converse assertion is obvious, using \reff{eq.prop.Pnx.recurrence.1}.
\qed

{\bf Remarks.}
1.  This proposition is also stated
by Viennot \cite[p.~III-18, Proposition~III.7]{Viennot_83},
where a combinatorial proof is sketched;
by Yang \cite[Theorem~2.3]{Yang_12};
by Cheon and Kim \cite[Theorem~4.1]{Cheon_15};
by Verde-Star \cite[Theorem~2.1]{Verde-Star_17};
and by Costabile, Gualtieri and Napoli \cite[Theorem~4.2]{Costabile_19}.

2. The recurrence \reff{eq.prop.Pnx.recurrence.2}
can also be written in vector form as $x \bfP = \Pi \bfP$,
where $\bfP$ is the column vector whose entries are the sequence
$(P_n(x))_{n \ge 0}$.
Iterating this, we see that
$x^r \, \bfP = \Pi^r \, \bfP$ for any integer $r \ge 0$,
or concretely
\be
   x^r P_n(x)   \;=\;  \sum_{j=0}^{n+r} (\Pi^r)_{nj} \, P_j(x)
   \;.
 \label{eq.prop.Pnx.recurrence.2.iterated}
\ee
It then follows that for any polynomial $q(x)$ we have
\be
   q(x) \, P_n(x)   \;=\;  \sum_{j=0}^{n+\deg q} (q(\Pi))_{nj} \, P_j(x)
   \;.
 \label{eq.prop.Pnx.recurrence.2.iterated2}
\vspace*{-4mm}
\ee
\myendremark

\medskip

The monic polynomials $P_n(x)$ defined by the recurrence
\reff{eq.prop.Pnx.recurrence.1} are also the characteristic polynomials
of the leading principal submatrices of the production matrix $\Pi$.
To state this result, let us introduce notation as follows:
For any matrix $A = (a_{ij})_{i,j \ge 0}$,
we write $A_n = (a_{ij})_{0 \le i,j \le n-1}$
for its $n \times n$ leading principal submatrix,
and $\Delta_n(A) = \det A_n$
for its $n \times n$ leading principal minor,
with the convention $\Delta_0(A) = 1$.
When the matrix is lower-Hessenberg,
these leading principal minors satisfy a recurrence
that is reasonably well known
\cite[pp.~251--252]{Franklin_68}
\cite[pp.~410--411]{Wilkinson_65},
though perhaps not as well known as it should be:

\begin{lemma}[Leading principal minors of a lower-Hessenberg matrix]
   \label{lemma.hessenberg}%
\hfill\break%
$\!\!$The leading principal minors $\Delta_n(H)$
of a lower-Hessenberg matrix $H = (h_{ij})_{i,j \ge 0}$
satisfy the recurrence
\be
   \Delta_n  \;=\;
   \sum_{j=0}^{n-1} (-1)^{n-1-j} \, h_{n-1,j} \,
                \Biggl( \prod_{i=j}^{n-2} h_{i,i+1} \!\Biggr) \,
                \Delta_j
   \;.
 \label{eq.lemma.hessenberg}
\ee
\end{lemma}

\proof
Laplace-expand $\det H_n$ in the last (i.e.\ $(n-1)$st) row.
When row $n-1$ and column $j$ are deleted from $H_n$, what remains is
\be
   \left[ \begin{array}{ccc|cccc}
                  & & & & & & \\
                  & \Scale[1.5]{H_j} &  & & & \\
                  & & & & & & \\[-2mm]
                  \hline
                  & & & & & & \\[-3mm]
                  \times & \cdots & \times & h_{j,j+1} & & &  \\
                  \times & \cdots & \times & \times & h_{j+1,j+2} & &  \\
                  \vdots &        & \vdots & \vdots & \vdots & \ddots & \\
                  \times & \cdots & \times & \times & \times & \cdots & h_{n-2,n-1} \\
          \end{array}
   \right] 
\ee
(blank entries are zero),
so that its determinant is
$\Bigl( \prod\limits_{i=j}^{n-2} h_{i,i+1} \!\Bigr) \, \Delta_j$.
\qed

\vspace*{-1mm}

For tridiagonal matrices, \reff{eq.lemma.hessenberg}
becomes a three-term recurrence that is much better known
\cite[p.~35]{Horn_13}.

\medskip

\begin{proposition}[Monic polynomials as characteristic polynomials of production matrix]
   \label{prop.charpoly}
Let $\Pi = (\pi_{nk})_{n,k \ge 0}$ be a unit-lower-Hessenberg matrix,
and let $(P_n(x))_{n \ge 0}$ be the sequence of monic polynomials
defined by the recurrence \reff{eq.prop.Pnx.recurrence.1}
with the initial condition $P_0(x) = 1$.
Then $P_n(x) = \det(xI - \Pi_n)$.
\end{proposition}

\proof
Applying \reff{eq.lemma.hessenberg} to the matrix $H = xI - \Pi$,
we obtain
\be
   \Delta_n  \;=\;  (x - \pi_{n-1,n-1}) \Delta_{n-1}
                    \:-\: \sum_{j=0}^{n-2} \pi_{n-1,j} \, \Delta_j
   \;,
\ee
which matches the recurrence \reff{eq.prop.Pnx.recurrence.1}.
\qed

{\bf Remarks.}
1. When the matrix $\Pi$ is tridiagonal, this result is classical
\cite[p.~26, Exercise~4.12]{Chihara_78}
\cite[p.~24, Theorem~2.2.4]{Ismail_05}.\footnote{
   I thank Alex Dyachenko for drawing my attention to this classical result,
   which inspired the generalization presented in
   Proposition~\ref{prop.charpoly}.
}

2. The general case of Proposition~\ref{prop.charpoly} is also known:
see, for instance,
\cite[eq.~(2.8) ff.]{Coussement_05},
\cite[Theorem~2.3]{Yang_12},
\cite[Corollary~3.1]{Verde-Star_17}
and \cite[Theorem~4.4]{Costabile_19};
see also \cite[Theorem~3.1]{Cheon_15}.

3. When $R = \R$ or $\C$, it follows from Proposition~\ref{prop.charpoly}
that the zeros of $P_n(x)$ are the eigenvalues of $\Pi_n$.
This suggests that the asymptotic zero distribution of the
polynomials $P_n(x)$ as $n \to\infty$ should be related to the
spectral properties of the infinite matrix $\Pi$
acting on a suitable space of sequences (for instance, $\ell^2(\N)$).
See e.g.\ 
\cite{Kalyagin_94,Kalyagin_95,Robert_03,Aptekarev_06,Zhang_11}.
\myendremark

\bigskip

Finally, the recurrence \reff{eq.prop.Pnx.recurrence.1} also leads to
a combinatorial formula, due to Viennot \cite[p.~III-16]{Viennot_83},
for the matrix elements of $B = \scro(\Pi)^{-1}$
in terms of those of $\Pi$:

\begin{corollary}[Viennot \protect\cite{Viennot_83}]
   \label{cor.Pnx.recurrence}
Let $\Pi = (\pi_{nk})_{n,k \ge 0}$ be a unit-lower-Hessenberg matrix,
and let $B = (b_{nj})_{n,j \ge 0}$ be given by $B = \scro(\Pi)^{-1}$.
Then $b_{nj}$ is the sum over partitions of $\{0,1,\ldots,n-1\}$
into zero or more intervals $[k,\ell]$ ($k \le \ell$)
and exactly $j$ empty sites,
with a weight $-\pi_{\ell k}$ for each interval $[k,\ell]$.\footnote{
   Viennot \cite[p.~III-16]{Viennot_83} inadvertently omitted
   the minus sign in front of $\pi_{\ell k}$ in this formula.
   But this was clearly an oversight, as he had the minus sign correct
   in the tridiagonal case \cite[p.~I-9]{Viennot_83}.
}
\end{corollary}


\proof
Write $\widehat{b}_{nj}$ for the quantity defined in the Corollary,
and define $\widehat{P}_n(x) = \sum_{j=0}^n \widehat{b}_{nj} \, x^j$.
Then $\widehat{P}_n(x)$ is the sum over partitions of $\{0,1,\ldots,n-1\}$
into zero or more intervals $[k,\ell]$ ($k \le \ell$)
and zero or more empty sites,
with a weight $-\pi_{\ell k}$ for each interval $[k,\ell]$
and a weight $x$ for each empty site.
And it is easy to see, by considering the status of the vertex $n$
in $\widehat{P}_{n+1}(x)$,
that the sequence $(\widehat{P}_n(x))_{n \ge 0}$
satisfies the same recurrence \reff{eq.prop.Pnx.recurrence.1}
as is satisfied by $(P_n(x))_{n \ge 0}$,
with the same initial condition $\widehat{P}_0(x) = P_0(x) = 1$.
So $\widehat{P}_n(x) = P_n(x)$.
\qed

%
%
%
\subsection{Summary}

To summarize the results obtained thus far:
There is a one-to-one correspondence between
sequences $(P_n(x))_{n \ge 0}$ of monic polynomials
(with coefficient matrix~$B$),
their dual sequences $(\scrl_k)_{k \ge 0}$ of linear functionals
(with moment matrix~$A$),
and their defining linear recurrences
\reff{eq.prop.Pnx.recurrence.1}/\reff{eq.prop.Pnx.recurrence.2}
(with production matrix~$\Pi$);
and these correspondences are given by \reff{eq.A.B.Pi}.

%
%
%
\subsection{Expectation values of products}

Fix now a unit-lower-Hessenberg matrix $\Pi$:
this defines (by Proposition~\ref{prop.Pnx.recurrence})
a sequence $(P_n(x))_{n \ge 0}$ of monic polynomials 
with coefficient matrix $B = \scro(\Pi)^{-1}$,
which in turn defines (by Proposition~\ref{prop.Pnx.dual})
a dual sequence $(\scrl_k)_{k \ge 0}$ of linear functionals 
with moment matrix $A = \scro(\Pi)$.
Our goal is to find a general combinatorial or algebraic formula
for quantities of the form $\scrl_k( q(x) \, P_m(x) \, P_n(x) )$,
where $q(x)$ is a polynomial,
in terms of the coefficients $\Pi$.
By analogy to probability theory,
we refer colloquially to quantities like
$\scrl_k( q(x) \, P_m(x) \, P_n(x) )$
as ``expectation values''.

In the tridiagonal case with $k = 0$,
Viennot \cite[p.~I-15, Proposition~I.17]{Viennot_83}
found a beautiful formula for these expectation values:

\begin{proposition}[Viennot \protect\cite{Viennot_83}]
   \label{prop.viennot}
When the unit-lower-Hessenberg matrix $\Pi$ is tridiagonal, we have
\be
   \scrl_0( q(x) \, P_m(x) \, P_n(x) )
   \;=\;
   \pi_{10} \pi_{21} \cdots \pi_{n,n-1} \, (q(\Pi))_{mn}
 \label{eq.prop.viennot}
\ee
for any polynomial $q(x)$.
In particular, when $q(x) = 1$ we have the orthogonality relation
\be
   \scrl_0( P_m(x) \, P_n(x) )
   \;=\;
   h_n \, \delta_{mn}
 \label{eq.prop.viennot.2}
\ee
with the normalizing constant $h_n = \pi_{10} \pi_{21} \cdots \pi_{n,n-1}$.
\end{proposition}

Please note \cite[p.~I-15]{Viennot_83}
that the right-hand side of \reff{eq.prop.viennot} with $q(x) = x^r$,
namely $\pi_{10} \pi_{21} \cdots \pi_{n,n-1} \, (\Pi^r)_{mn}$,
can be interpreted as the total weight for Motzkin paths of length $r+m+n$
from height $0 \to 0$ in which the first $m$ steps are ``up'' steps
(getting weight~1) and the last $n$ steps are ``down'' steps
(getting weight $\pi_{n,n-1} \cdots \pi_{10}$).
Now, in a Motzkin path that starts and ends at the same height,
each ``up'' step $i \to i+1$ can be paired with a ``down'' step $i+1 \to i$;
it~follows that the weight of such a path
equals the weight of the reversed path.
These considerations show that the right-hand side of \reff{eq.prop.viennot}
is indeed symmetric in $m \leftrightarrow n$.

Viennot \cite[pp.~I-16--I-19]{Viennot_83}
proved Proposition~\ref{prop.viennot}
by a rather intricate combinatorial argument;
here we give a simple algebraic proof:

\proofof{Proposition~\ref{prop.viennot}}
Let $\bll = (\ell_n)_{n \ge 0}$ be the moment sequence of
the linear functional $\scrl_0$, i.e.
\be
   \ell_n  \;=\;  \scrl_0(x^n)  \;=\;  a_{n0}  \;=\; (\Pi^n)_{00}
   \;.
\ee
And let $H_\infty(\bll) = (\ell_{i+j})_{i,j \ge 0}$
be the Hankel matrix associated to the sequence $\bll$.
Since the unit-lower-Hessenberg matrix $\Pi$ is tridiagonal,
it is of the form \reff{def.Jnl.production.bis.0}
with $\gamma_n = \pi_{nn}$ and $\beta_n = \pi_{n,n-1}$.
Proposition~\ref{prop.LDLT} therefore gives
\be
   H_\infty(\bll)  \;=\; A D A^{\rm T}
\ee
where $D = \diag(1,\beta_1, \beta_1 \beta_2, \ldots)$.
Since $B = A^{-1}$, we can rewrite this as
\be
   B \, H_\infty(\bll) \, B^{\rm T}  \;=\;  D
   \;,
 \label{eq.proof.prop.viennot.1}
\ee
or concretely
\be
   \sum_{i,j=0}^\infty  b_{mi} \, b_{nj} \, \ell_{i+j}
   \;=\;
   \beta_1 \cdots \beta_n \, \delta_{mn}
   \;.
 \label{eq.proof.prop.viennot.2}
\ee
But the left-hand side of \reff{eq.proof.prop.viennot.2}
is exactly $\scrl_0( P_m(x) \, P_n(x) )$.
This proves \reff{eq.prop.viennot.2}.

If we further left-multiply \reff{eq.proof.prop.viennot.1} by $q(\Pi)$,
we obtain
\be
   q(\Pi) \, B \, H_\infty(\bll) \, B^{\rm T}  \;=\;  q(\Pi) \, D
   \;,
 \label{eq.proof.prop.viennot.3}
\ee
or concretely
\be
   \sum_{i,j,m'=0}^\infty  (q(\Pi))_{mm'} \, b_{m' i} \, b_{nj} \, \ell_{i+j}
   \;=\;
   (q(\Pi))_{mn} \, \beta_1 \cdots \beta_n
   \;.
 \label{eq.proof.prop.viennot.4}
\ee
But, using \reff{eq.prop.Pnx.recurrence.2.iterated2},
we see that the left-hand side of \reff{eq.proof.prop.viennot.4}
is exactly $\scrl_0( q(x) \, P_m(x) \, P_n(x) )$.
This proves \reff{eq.prop.viennot}.
\qed

H\'elder Lima \cite{Lima_private}
has pointed out to me that Proposition~\ref{prop.viennot}
can be extended to $k \neq 0$ as follows:

\begin{proposition}
   \label{prop.viennot2}
When the unit-lower-Hessenberg matrix $\Pi$ is tridiagonal, we have
\be
   \scrl_k( q(x) )
   \;=\;
   h_k^{-1} \, \scrl_0( q(x) \, P_k(x))
 \label{eq.prop.viennot2.1}
\ee
and hence
\be
   \scrl_k( q(x) \, P_m(x) \, P_n(x) )
   \;=\;
   h_k^{-1} h_n \, (q(\Pi) \, P_k(\Pi))_{mn}
 \label{eq.prop.viennot2.2}
\ee
for any polynomial $q(x)$,
where $h_n = \pi_{10} \pi_{21} \cdots \pi_{n,n-1}$.
\end{proposition}

\proof
Put $m=0$ in \reff{eq.prop.viennot}, rename $n$ as $k$, and take $q(x) = x^n$:
this gives
\be
   \scrl_0(x^n \, P_k(x))
   \;=\;
   h_k \, (\Pi^n)_{0k}
   \;=\;
   h_k \, a_{nk}
   \;=\;
   h_k \, \scrl_k(x^n)
   \;.
\ee
It follows that
\be
   \scrl_0(q(x) \, P_k(x))
   \;=\;
   h_k \, \scrl_k(q(x))
\ee
for any polynomial $q(x)$.
This proves \reff{eq.prop.viennot2.1}.
Then replace $q(x)$ by $q(x) \, P_m(x) \, P_n(x)$
and use \reff{eq.prop.viennot};
this proves \reff{eq.prop.viennot2.2}.
\qed

{\bf Remark.}
The meaning of $h_k^{-1}$ in the formulae
\reff{eq.prop.viennot2.1} and \reff{eq.prop.viennot2.2}
requires some clarification.
The formulae obviously hold if $h_k$ is invertible in the ring $R$.
In particular this is the case if we consider
$\bm{\Pi} = \{\pi_{ij}\}_{i \ge j \ge 0}$ to be indeterminates
and we work in the ring $\Z[\bm{\Pi},\bm{\Pi}^{-1}]$ of Laurent polynomials.
Now suppose that $k \le n$:
then both sides of \reff{eq.prop.viennot2.2}
are in fact polynomials in $\bm{\Pi}$,
so the identity holds when the $\pi_{ij}$ are specialized
to arbitrary elements in an arbitrary commutative ring.
I~am not sure what happens when $k > n$.
\myendremark

\medskip

We would now like to generalize these results to the non-tridiagonal case.
(Viennot \cite[top p.~V-41]{Viennot_83} alludes to this as an open problem.)
For the moment I have only the following result:

\begin{theorem}
   \label{thm.gen_viennot}
With $(P_n(x))_{n \ge 0}$ and $(\scrl_k)_{k \ge 0}$ defined as above,
we have
\be
   \scrl_k( x^\ell \, P_n(x) )
   \;=\;
   (\Pi^\ell)_{nk}
      \label{eq.thm.gen_viennot.0}
\ee
and more generally
\begin{subeqnarray}
   \scrl_k( x^\ell \, P_m(x) \, P_n(x) )
   & = &
   \sum_{j=0}^m b_{mj} \, (\Pi^{j+\ell})_{nk}
      \slabel{eq.thm.gen_viennot.1} \\[2mm]
   & = &
   (B H_{k}^{(\ell)} B^{\rm T})_{mn}
      \slabel{eq.thm.gen_viennot.2}
\end{subeqnarray}
where $H_{k}^{(\ell)}$ is the $\ell$-shifted Hankel matrix
of the $k^{th}$ column of $A$:
\be
   (H_{k}^{(\ell)})_{ij}  \;\eqdef\;  a_{i+j+\ell,k}  \;=\;  (\Pi^{i+j+\ell})_{0k}
   \;.
 \label{def.Hkl}
\ee
\end{theorem}

\proof
Since $(P_n(x))_{n \ge 0}$ and $(\scrl_k)_{k \ge 0}$ are a dual pair,
we have $\scrl_k(P_n(x)) = \delta_{kn}$.
Applying $\scrl_k$ to \reff{eq.prop.Pnx.recurrence.2.iterated},
we obtain $\scrl_k(x^\ell P_n(x)) = (\Pi^\ell)_{nk}$;
this proves \reff{eq.thm.gen_viennot.0}.
Now use $P_m(x) = \sum_{j=0}^m b_{mj} \, x^j$;
inserting this into \reff{eq.thm.gen_viennot.0}
gives \reff{eq.thm.gen_viennot.1}.

On the other hand, we have $\scrl_k(x^n) = a_{nk}$.
Inserting the representing equation in terms of $B$
for both $P_m(x)$ and $P_n(x)$, we have
\begin{subeqnarray}
   \scrl_k( x^\ell \, P_m(x) \, P_n(x) )
   & = &
   \scrl_k\biggl( \sum\limits_{i=0}^m \sum\limits_{j=0}^n
                   b_{mi} \, b_{nj} \, x^{i+j+\ell}
          \biggr)
      \\[2mm]
   & = &
   \sum\limits_{i=0}^m \sum\limits_{j=0}^n
                   b_{mi} \, b_{nj} \, a_{i+j+\ell,k}
      \\[2mm]
   & = &
   (B H_{k}^{(\ell)} B^{\rm T})_{mn}
\end{subeqnarray}
when $H_{k}^{(\ell)}$ is defined by \reff{def.Hkl}.
This proves \reff{eq.thm.gen_viennot.2}.
\qed

{\bf Remark.}  The identity \reff{eq.thm.gen_viennot.0}
is contained in the thesis of Roblet \cite[p.~153, Proposition~78]{Roblet_94}
in the special case where $\Pi$ is $(d,1)$-banded and $0 \le k \le d-1$
(but this is no real loss of generality, since we can take $d \to\infty$).
Roblet's proof was combinatorial,
following the model of Viennot \cite[pp.~I-15--I.19]{Viennot_83}.
\myendremark

It is an immediate consequence of Theorem~\ref{thm.gen_viennot}
that certain matrix elements have to vanish:

\begin{corollary}
   \label{cor.gen_viennot}
In the situation of Theorem~\ref{thm.gen_viennot}:
\begin{itemize}
   \item[(a)]  $\scrl_k( x^\ell \, P_m(x) \, P_n(x)) \,=\, 0$
       whenever $k > \ell+m+n$.
   \item[(b)]  If $\Pi$ is $(d,1)$-banded
       [that is, $\pi_{nk} = 0$ whenever $k < n-d$],
       then $\scrl_k( x^\ell \, P_m(x) \, P_n(x)) \,=\, 0$
       whenever $k < n - d(\ell+m)$ or $k < m - d(\ell+n)$.
\end{itemize}
Conversely, if $\scrl_k( x P_n(x)) \,=\, 0$ whenever $k < n-d$
[this is the $\ell=1$, $m=0$ case of (b)],
then $\Pi$ is $(d,1)$-banded.
\end{corollary}

\proof
(a) is an immediate consequence of \reff{eq.thm.gen_viennot.1}
together with the facts that $B$ is lower-triangular
and $\Pi$ is lower-Hessenberg.

(b) The vanishing for $k < n - d(\ell+m)$ is likewise
an immediate consequence of \reff{eq.thm.gen_viennot.1}
together with the facts that $B$ is lower-triangular
and $\Pi$ is $(d,1)$-banded.
The vanishing for $k < m - d(\ell+n)$ then follows
from the symmetry $m \leftrightarrow n$
of the left-hand side of \reff{eq.thm.gen_viennot.1}.

The converse assertion follows trivially from \reff{eq.thm.gen_viennot.0}.
\qed

In particular, if $\Pi$ is $(d,1)$-banded,
then $\scrl_k( x^\ell \, P_n(x)) = 0$ whenever $n > d\ell+k$.
So the sequence $(P_n(x))_{n \ge 0}$ is \textbfit{$\bm{d}$-orthogonal}
with respect to the sequence $\scrl_0,\ldots,\scrl_{d-1}$ of linear forms,
in the sense defined in the Introduction.
But in fact this vanishing for $n > d\ell+k$ holds
for {\em all}\/ $k \ge 0$, not just for $k \le d-1$.

\bigskip

Unfortunately equation~\reff{eq.thm.gen_viennot.1} is not very nice,
because it is not manifestly symmetric in $m \leftrightarrow n$;
moreover, it is not expressed solely in terms of $\Pi$.
The equation~\reff{eq.thm.gen_viennot.2} has the $m \leftrightarrow n$ symmetry,
but it is still not expressed solely in terms of $\Pi$.
Combining \reff{eq.thm.gen_viennot.1} with Corollary~\ref{cor.Pnx.recurrence}
gives an explicit formula for $\scrl_k( x^\ell \, P_m(x) \, P_n(x) )$
in terms of the matrix elements of $\Pi$;
but this formula is rather complicated,
and it also fails to make manifest the symmetry $m \leftrightarrow n$.
Combining \reff{eq.thm.gen_viennot.2} and \reff{def.Hkl}
with Corollary~\ref{cor.Pnx.recurrence}
gives an even more complicated formula
for $\scrl_k( x^\ell \, P_m(x) \, P_n(x) )$
in terms of the matrix elements of $\Pi$,
which at least is manifestly symmetric in $m \leftrightarrow n$.
But none of these formulae seem really satisfactory.
What we really want is something that looks more like
Propositions~\ref{prop.viennot} and \ref{prop.viennot2}
and that reduces to them when $\Pi$ is tridiagonal.
My hope is that there might be some cancellations between terms
in the expansion of \reff{eq.thm.gen_viennot.1},
so that the final result can perhaps be stated in a simpler way.
I~therefore conclude by stating the main unsolved problem of this paper:

\begin{openproblem}
   \label{openproblem.expectation}
Find a more satisfactory formula for $\scrl_k( x^\ell \, P_m(x) \, P_n(x) )$
--- ideally one that resembles
Propositions~\ref{prop.viennot} and \ref{prop.viennot2}
and that reduces to them when $\Pi$ is tridiagonal.
Even the special case $k=\ell=0$ would be of great interest.
\end{openproblem}

\section{Sequence of monic polynomials orthogonal to a sequence of
           linear functionals}   \label{sec.orthogonal}

Let $\bGamma = (\Gamma_k)_{k \ge 0}$
be a sequence of linear functionals on $R[x]$,
with moment matrix $\Gamma = (\gamma_{nk})_{n,k \ge 0}$
given by $\gamma_{nk} = \Gamma_k(x^n)$.
And let $\bfP = (P_n(x))_{n \ge 0}$ be a sequence of monic polynomials,
with unit-lower-triangular coefficient matrix $B = (b_{nj})_{n,j \ge 0}$
given by $P_n(x) = \sum_{j=0}^n b_{nj} \, x^j$.
We say that \textbfit{$\bfP$ is orthogonal to $\bGamma$}
in case each $P_n(x)$ is orthogonal to $\Gamma_k$ for $0 \le k \le n-1$,
i.e.\ $\Gamma_k(P_n(x)) = 0$ for $0 \le k \le n-1$.
Since
\be
   \Gamma_k(P_n(x))  \;=\;  \sum_{j=0}^n b_{nj} \, \gamma_{jk}
                        \;=\;  (B\Gamma)_{nk}
   \;,
\ee
we see that $\bfP$ is orthogonal to $\bGamma$
if and only if $B\Gamma$ vanishes below the diagonal,
or in other words $B\Gamma$ is an upper-triangular matrix $U$,
or equivalently $\Gamma = B^{-1} U$.
We record this simple fact:

\begin{proposition}[Orthogonality between a sequence of linear functionals and
   a sequence of monic polynomials]
   \label{prop.ortho}
Let $\bGamma = (\Gamma_k)_{k \ge 0}$
be a sequence of linear functionals on $R[x]$,
with moment matrix $\Gamma$;
and let $\bfP = (P_n(x))_{n \ge 0}$ be a sequence of monic polynomials,
with unit-lower-triangular coefficient matrix $B$.
Then $\bfP$ is orthogonal to $\bGamma$
if and only if the matrix $B\Gamma$ is upper-triangular,
or equivalently if there exists an upper-triangular matrix $U$
such that $\Gamma = B^{-1} U$.
\end{proposition}

We can interpret this result in two ways,
depending on whether we start from $\bGamma$ or from $\bfP$:

\medskip

{\bf Starting from $\bGamma$.}
Suppose that the moment matrix $\Gamma$
has a factorization $\Gamma = LU$
where $L$ is unit-lower-triangular and $U$ is upper-triangular.
Then $B = L^{-1}$ is the coefficient matrix
for a sequence $\bfP$ of monic polynomials that is orthogonal to $\bGamma$.
Furthermore, the sequence $\bfL = (\scrl_k)_{k \ge 0}$
of linear functionals dual to $\bfP$
(given by Proposition~\ref{prop.Pnx.dual})
then has moment matrix $A = B^{-1} = L$.
It follows that $\Gamma_k = \sum_{j=0}^k (U^{\rm T})_{kj} \, \scrl_j$;
and if the diagonal elements of $U$ are invertible,
then $\scrl_k = \sum_{j=0}^k (U^{-\rm T})_{kj} \, \Gamma_j$.
In other words, each $\Gamma_k$ is a linear combination of
$\scrl_0,\ldots,\scrl_k$;
and if the diagonal elements of $U$ are invertible,
then each $\scrl_k$ is a special linear combination
of $\Gamma_0,\ldots,\Gamma_k$, namely,
one that makes the first $k$ elements of its moment sequence zero
and the next element 1.

\medskip

{\bf Starting from $\bfP$.}
Let $\bfP = (P_n(x))_{n \ge 0}$ be a sequence of monic polynomials,
with coefficient matrix $B$.
Then there is a {\em canonically associated}\/ sequence of linear functionals
with respect to which $\bfP$ is orthogonal,
namely, the dual sequence $\bfL = (\scrl_k)_{k \ge 0}$,
with moment matrix $A = B^{-1}$.
But $\bfP$ is also orthogonal with respect to {\em any}\/ sequence
$\bGamma = (\Gamma_k)_{k \ge 0}$ of linear functionals
whose moment matrix $\Gamma$ is of the form $\Gamma = AU$,
where $U$ is {\em any}\/ upper-triangular matrix;
or in other words, each $\Gamma_k$ is an arbitrary linear combination
of $\scrl_0,\ldots,\scrl_k$.

\medskip

{\bf Remark.}  When $R = \R$, it is often the case in applications
that $\scrl_0$ is the moment functional of a positive measure.
But $\scrl_k$ for $k \ge 1$ {\em cannot}\/ be the moment functional
of a positive measure, because the zeroth component
of its moment sequence vanishes
but the sequence is not identically zero (since the $k$th component is 1).
Analyst readers may be interested in the following problem:
Given a unit-lower-triangular matrix $A$
whose zeroth column is the moment sequence of a positive measure,
find an upper-triangular matrix $U$
such that every column of $\Gamma = AU$
is the moment sequence of a positive measure;
or in other words, find linear combinations $\Gamma_k$
of $\scrl_0,\ldots,\scrl_k$ that are all moment functionals
of positive measures.
This problem of course has trivial solutions:
we could take all the $\Gamma_k$ to be zero, or to be equal to $\scrl_0$.
But suppose we further insist that
all the functionals $\Gamma_0,\Gamma_1,\ldots$ be linearly independent.
Then the problem seems to be nontrivial.
\myendremark

\medskip

Since $\Gamma_k$ is a linear combination of $\scrl_0,\ldots,\scrl_k$
whenever $\bfP$ is orthogonal to $\bGamma$,
it also follows that the result of Corollary~\ref{cor.gen_viennot}(b)
holds with $\scrl_k$ replaced by $\Gamma_k$:

\begin{corollary}
   \label{cor.gen_viennot.bis}
Let $\Pi = (\pi_{nk})_{n,k \ge 0}$ be a unit-lower-Hessenberg matrix,
let $\bfP = (P_n(x))_{n \ge 0}$ be the sequence of monic polynomials
defined by \reff{eq.prop.Pnx.recurrence.1}/\reff{eq.prop.Pnx.recurrence.2}
with initial condition $P_0(x) = 1$,
and let $\bGamma = (\Gamma_k)_{k \ge 0}$ be any sequence of linear functionals
such that $\bfP$ is orthogonal to $\bGamma$.

If $\Pi$ is $(d,1)$-banded,
then $\Gamma_k( x^\ell \, P_m(x) \, P_n(x)) \,=\, 0$
whenever $k < n - d(\ell+m)$ or $k < m - d(\ell+n)$.
\end{corollary}

In view of the key role played here by the factorization $\Gamma = LU$
when we start from $\bGamma$,
it is now appropriate to recall some simple facts
concerning the existence and uniqueness of $LU$ factorizations
for matrices over a commutative ring \cite{Sokal_LU}.
Let us say that a square matrix $\Gamma$ has a
\textbfit{weak $\bm{LU}$ factorization} if $\Gamma = LU$ where
$L$ is lower-triangular and $U$ is upper-triangular,
and an \textbfit{$\bm{LU}$ factorization} if $\Gamma = LU$ where
$L$ is unit-lower-triangular and $U$ is upper-triangular.
We then have \cite{Sokal_LU}:

\begin{proposition}[$LU$ factorization for matrices over a commutative ring]
   \label{prop.LU}
Let $\Gamma$ be an $n \times n$ matrix with entries in a commutative ring $R$,
and let $\Delta_1,\ldots,\Delta_n$ be its leading principal minors.
\begin{itemize}
   \item[(a)]  If $\Gamma$ has a weak $LU$ factorization,
       then we must have
       $\Delta_1 \divides \Delta_2 \divides \cdots \divides \Delta_n$,
       where $a \divides b$ denotes that $a$~\hbox{divides}~$b$
       in the ring $R$.
   \item[(b)]  If $\Gamma$ has a weak $LU$ factorization
       in which none of the diagonal elements of $L$ or $U$
       is a zero or a divisor of zero, 
       then none of $\Delta_1,\ldots,\Delta_n$
       is zero or a divisor of zero.
   \item[(c)]  If $\Gamma$ has a weak $LU$ factorization
       in which all of the diagonal elements of $L$ and $U$
       are invertible in $R$,
       then $\Delta_1,\ldots,\Delta_n$ are invertible in $R$.
\end{itemize}
Conversely,
\begin{itemize}
   \item[(d)]  If none of $\Delta_1,\ldots,\Delta_{n-1}$
       is a zero or a divisor of zero,
       then $\Gamma$ has \emph{at~most~one} $LU$ factorization.
       (In particular this holds if $R$ is an integral domain
        and $\Delta_1,\ldots,\Delta_{n-1} \neq 0$.)
   \item[(e)]  If $\Delta_1,\ldots,\Delta_{n-1}$ are invertible in $R$,
       then $\Gamma$ has \emph{exactly~one} $LU$ factorization.
       (In particular this holds if $R$ is a field
        and $\Delta_1,\ldots,\Delta_{n-1} \neq 0$.)
\end{itemize}
\end{proposition}

\noindent
Since \cite{Sokal_LU} is not yet publicly available,
we include a proof of Proposition~\ref{prop.LU} in the Appendix.

Taking $n \to\infty$ in Proposition~\ref{prop.LU}(d,e)
and applying it to the situation considered
in Proposition~\ref{prop.ortho}, we conclude:

\begin{corollary}[Existence and uniqueness of a sequence of monic
   polynomials\break orthogonal to a given sequence of linear functionals]
   \label{cor.LU}
Let $R$ be a commutative ring, and let $\bGamma = (\Gamma_k)_{k \ge 0}$
be a sequence of linear functionals on $R[x]$,
with moment matrix $\Gamma$.
\begin{itemize}
   \item[(a)] If none of the leading principal minors
$\Delta_1,\Delta_2,\ldots$ of $\Gamma$ is zero or a divisor of zero,
then there is \emph{at~most~one} sequence of monic polynomials
orthogonal to $\bGamma$.
(In particular this holds if $R$ is an integral domain
 and $\Delta_1,\Delta_2,\ldots \neq 0$.)
   \item[(b)] If all of the leading principal minors
$\Delta_1,\Delta_2,\ldots$ of $\Gamma$ are invertible in $R$,
then there is \emph{exactly~one} sequence of monic polynomials
orthogonal to $\bGamma$.
(In particular this holds if $R$ is a field
 and $\Delta_1,\Delta_2,\ldots \neq 0$.)
\end{itemize}
\end{corollary}

\section{Application to ordinary orthogonal polynomials}
   \label{sec.prodmat_OOP}

Let us begin by showing how the general theory from the preceding section
applies to ordinary orthogonal polynomials.

Fix a linear functional $\scrl$,
with moment sequence $\bll = (\ell_n)_{n \ge 0}$
given by $\ell_n = \scrl(x^n)$.
And let us choose $\Gamma_k$ to be the $k$-shift of $\scrl$:
that is, $\Gamma_k(x^n) \eqdef \scrl(x^{n+k})$.
Then the moment matrix $\Gamma$ 
of the sequence $\bGamma = (\Gamma_k)_{k \ge 0}$ of linear functionals
is the Hankel matrix $H_\infty(\bll) = (\ell_{i+j})_{i,j \ge 0}$
associated to the sequence $\bll$:
that is, $\gamma_{nk} = \ell_{n+k}$.
And a sequence $\bfP = (P_n(x))_{n \ge 0}$ of monic polynomials
is orthogonal to $\bGamma$ in case $\scrl(x^k \, P_n(x)) = 0$
for $0 \le k \le n-1$,
i.e.\ precisely when $\bfP$ is a sequence of monic orthogonal polynomials
in the usual sense associated to the linear functional $\scrl$
\cite[Chapter~1]{Chihara_78}.
In particular, by Corollary~\ref{cor.LU}(b),
such a sequence $\bfP$ exists (and is unique)
whenever $R$ is a field and
all the leading principal minors $\Delta_1,\Delta_2,\ldots$
of $\Gamma$ are nonzero.

Let us now relate this to production matrices
and classical continued fractions.
It is known
\cite[Theorem~51.1]{Wall_48}
\cite[p.~IV-17, Corollaire~7 and p.~V-5, Proposition~1]{Viennot_83}
that if $R$ is a field and
all the leading principal minors $\Delta_1,\Delta_2,\ldots$
of the Hankel matrix $\Gamma = H_\infty(\bll)$ are nonzero,
then there exists a classical J-fraction that represents
the ordinary generating function of the sequence $\bll$
that underlies this Hankel matrix, i.e.
\be
   \sum_{n=0}^\infty \ell_n \, t^n
   \;=\;
   \cfrac{1}{1 - \gamma_0 t - \cfrac{\beta_1 t^2}{1 - \gamma_1 t - \cfrac{\beta_2 t^2}{1 - \gamma_2 t - \cfrac{\beta_3 t^2}{1- \gamma_3 t - \cdots}}}}
 \label{eq.Jtype.cfrac}
\ee
in the sense of formal power series,
with coefficients $\gamma_0,\gamma_1,\ldots \in R$
and $\beta_1,\beta_2,\ldots \in R \setminus \{0\}$.
In fact, the J-fraction coefficients are connected to the
leading principal minors by
\be
   \beta_n  \;=\;  {\Delta_{n-1} \, \Delta_{n+1} \over \Delta_n^2}
\ee
where $\Delta_{-1} \eqdef 1$ [this follows from \reff{eq.deltas.Jfrac.0}],
while the $\gamma_n$ are given by other determinants
involving the moments $\bll$
\cite[Sections~IV.3 and V.1]{Viennot_83}.\footnote{
   Please note that Viennot's $\Delta_n$ \cite[p.~IV-15, eqn.~(17)]{Viennot_83}
   is an $(n+1) \times (n+1)$ determinant,
   hence equal to my $\Delta_{n+1}$.
}

This J-fraction has all the properties described
in Section~\ref{subsec.classical.CF}.
In particular, it has a tridiagonal production matrix $\Pi$
in which $\pi_{n,n+1} = 1$, $\pi_{nn} = \gamma_n$, $\pi_{n,n-1} = \beta_n$
and $\pi_{nk} = 0$ for $k < n-1$ or $k > n+1$
[cf.\ \reff{def.Jnl.production.bis.0}].
The zeroth column of the output matrix $\sfJ = \scro(\Pi)$
is the moment sequence $\bll$.
Furthermore, Proposition~\ref{prop.LDLT} tells us that
the Hankel matrix $\Gamma = H_\infty(\bll)$ has the $LDL^{\rm T}$ factorization
\be
   \Gamma  \;=\;  \sfJ D \sfJ^{\rm T}
   \;,
  \label{eq.jacobi.LDLT.bis}
\ee
where $\sfJ = \scro(\Pi)$
is the unit-lower-triangular matrix of generalized Jacobi--Rogers polynomials,
and $D = \diag(1,\beta_1, \beta_1 \beta_2, \ldots)$.
By Proposition~\ref{prop.LU}(d,e),
this gives the {\em unique}\/ $LU$ factorization of $\Gamma$,
i.e.\ $L = \sfJ$ and $U = D \sfJ^{\rm T}$.
Proposition~\ref{prop.ortho} then implies that there is
a unique sequence $\bfP$ of monic polynomials orthogonal to $\bGamma$,
and its coefficient matrix is $B = \sfJ^{-1} = \scro(\Pi)^{-1}$.
So the $LDL^{\rm T}$ factorization can be written as
\be
   \Gamma  \;=\;  A D A^{\rm T}  \;=\;  B^{-1} D B^{-\rm T}
\ee
where $A = \scro(\Pi)$ and $B = A^{-1} = \scro(\Pi)^{-1}$.
(This is the general factorization $\Gamma = B^{-1} U$
found in Proposition~\ref{prop.ortho},
specialized to a case in which the matrix $\Gamma$ is symmetric.)

Finally, Proposition~\ref{prop.Pnx.recurrence}
implies that the orthogonal polynomials obey the three-term recurrence
\be
   P_{n+1}(x)
   \;=\;
   (x - \gamma_n) \, P_n(x) \:-\: \beta_n \, P_{n-1}(x)
   \;,
\ee
where the coefficients arising in the recurrence are precisely the same ones
that arise in the J-fraction for the moment sequence $\bll$.
This is, of course, a well-known fact
\cite[Theorems~50.1 and 51.1]{Wall_48}
\cite[Chapitre~V]{Viennot_83}
\cite[Theorem~2.3 and Corollary~2.5]{Corteel_16}
\cite[Section~5.2.1]{Zeng_21}.
And it is likewise well known that
the coefficient matrix of the orthogonal polynomials
is $B = A^{-1} = \scro(\Pi)^{-1}$
\cite[p.~III-2, Th\'eor\`eme~1]{Viennot_83} \cite[Proposition~5.12]{Zeng_21}.
But it is pleasing to see all these classical facts come together
as consequences of a simple algebraic theory.

\section{Application to multiple orthogonal polynomials}
   \label{sec.prodmat_MOP}

Let us now apply the general theory from the Section~\ref{sec.orthogonal}
to the multiple orthogonal polynomials of type~II
along an increasing nearest-neighbor path in $\N^r$.

If $\scrl$ is a linear functional on $R[x]$
and $k$ is a nonnegative integer,
we denote by $\scrl^{\sharp k}$ the $k$-shift of $\scrl$:
that is, $\scrl^{\sharp k}(x^n) \eqdef \scrl(x^{n+k})$.

Now fix an integer $r \ge 1$,
and fix linear functionals $\scrl^{(1)},\ldots,\scrl^{(r)}$ on $R[x]$.
(In the analytical setting, we will have $R = \R$,
 and $\scrl^{(1)},\ldots,\scrl^{(r)}$ will be the moment functionals
 associated to the positive measures $\mu_1,\ldots,\mu_r$.)
And let $(P_{\bfn}(x))_{\bfn \in \N^r}$ be the
multiple orthogonal polynomials of type~II
associated to the linear functionals $\scrl^{(1)},\ldots,\scrl^{(r)}$,
which we here {\em assume}\/ to exist.

Now let $j_1,j_2,\ldots$ be an infinite sequence
of elements of $\{1,\ldots,r\}$,
and define a sequence $(\bfn_k)_{k \ge 0}$ of multi-indices in $\N^r$
by $\bfn_k = \sum_{i=1}^k \bfe_{j_i}$.
They satisfy $|\bfn_k| = k$
and describe an increasing nearest-neighbor path in $\N^r$
in which the $i$th step is along direction $j_i$.
Now let $\Phat_k(x) \eqdef P_{\bfn_k}(x)$
be the multiple orthogonal polynomial of type~II along this path in $\N^r$.
Let $m_i = (\bfn_{i-1})_{j_i} = (\bfn_i)_{j_i} - 1$
be the number of indices $j_1,\ldots,j_{i-1}$ that equal $j_i$.
Then $\Phat_k(x)$ is orthogonal to the linear functionals
$\scrl^{\star 1},\ldots,\scrl^{\star k}$,
where the ``new'' linear functional appearing at stage $k$ is
\be
   \scrl^{\star k}
   \;=\;
   (\scrl^{(j_k)})^{\sharp\, m_k}
   \;,
\ee
i.e.
\be
   \scrl^{\star k}(x^n)
   \;=\;
   \scrl^{(j_k)}(x^{n + m_k})
   \;.
\ee
Now set $\Gamma_k = \scrl^{\star,k+1}$:
we then see that the sequence $\bfPhat = (\Phat_k(x))_{k \ge 0}$
is orthogonal to the sequence $\bGamma = (\Gamma_k)_{k \ge 0}$
in the sense of the preceding section.

On the other hand, we know from the general theory
of multiple orthogonal polynomials (Section~\ref{subsec.MOP})
that the sequence $\bfPhat$ satisfies an $(r+2)$-term linear recurrence
\reff{eq.prop.Ptildenx.recurrence.2}
with an $(r,1)$-banded unit-lower-Hessenberg matrix $\Pi$.
It follows from Proposition~\ref{prop.Pnx.recurrence}
that the coefficient matrix of the sequence $\bfPhat$
is the unit-lower-triangular matrix $B = \scro(\Pi)^{-1}$.
Proposition~\ref{prop.ortho} then implies that
$\Gamma = B^{-1} U = \scro(\Pi) \, U$
for some upper-triangular matrix $U$;
this is the $LU$ factorization of $\Gamma$.
Thus, $\Gamma_0$ is proportional to the zeroth column of
the output matrix $\scro(\Pi)$;
and more generally, each $\Gamma_k$ is a linear combination of
columns $0,\ldots,k$ of $\scro(\Pi$).

\begin{example}[Multiple orthogonal polynomials along an axis]
   \label{exam.MOP_axis}
\rm
If $j_1,j_2,j_3,\ldots = j,j,j,\ldots\,$,
then $\bfn_k = k \bfe_j$,
and we are in the situation of Section~\ref{sec.prodmat_OOP}
with $\scrl = \scrl^{(j)}$ and $\Gamma_k = (\scrl^{(j)})^{\sharp k}$.
Then $(\Phat_k(x))_{k \ge 0}$
is the sequence of ordinary orthogonal polynomials
associated to the linear functional $\scrl^{(j)}$.
\myendremark
\end{example}

\begin{example}[Multiple orthogonal polynomials along the stepline]
   \label{exam.MOP_stepline}
\rm
If $j_1,j_2,j_3,\ldots = 1,\ldots,r,1,\ldots,r,\ldots\,$,
then $(\bfn_k)_{k \ge 0}$ is the stepline
defined in \reff{def.Ptilde.stepline},
i.e.\ $\Phat_k(x) = \Ptilde_k(x)$.
The sequence $\scrl^{\star 1}, \scrl^{\star 2}, \ldots$
is $\scrl^{(1)},\ldots,\scrl^{(r)}$,
$(\scrl^{(1)})^{\sharp 1},\ldots,(\scrl^{(r)})^{\sharp 1}$,
$(\scrl^{(1)})^{\sharp 2},\ldots,(\scrl^{(r)})^{\sharp 2}$, \ldots\ .


In the case $r=2$, the formulae for $\scrl_0$ and $\scrl_1$
(but not the higher $\scrl_k$) can be found already in the thesis of Drake
\cite[Theorem~2.4.1]{Drake_06}.
\myendremark
\end{example}

\section{Some examples}  \label{sec.examples}

\subsection[Bessel $K_\nu$ weights $\Rightarrow$ Rising-factorial moments]{Bessel $\bm{K_\nu}$ weights $\Rightarrow$ Rising-factorial moments}
   \label{subsec.BesselK}

Two decades ago, {Van Assche} and Yakubovich \cite{VanAssche_00}
studied the multiple orthogonal polynomials of types~I and II
associated to a pair of measures (that is, $r=2$)
in which the weights are modified Bessel functions of the second kind,
$K_\nu(x)$ \cite[p.~78]{Watson_44}, multiplied by powers of $x$.
We shall follow their paper closely,
but change the notation to make the formulae more symmetrical.

For real numbers $a_1, a_2 > 0$,
let $\mu_{a_1,a_2}$ be the positive measure on $[0,\infty)$ given by
\be
   d\mu_{a_1,a_2}(x)
   \;=\;
   {2 \over \Gamma(a_1) \, \Gamma(a_2)} \:
   x^{(a_1 + a_2 - 2)/2} \, K_{a_1 - a_2}(2\sqrt{x}) \: dx
   \;.
 \label{def.besselK.weight}
\ee
This is symmetric in $a_1 \leftrightarrow a_2$
because $K_{-\nu} = K_\nu$.
The moments of the measure $\mu_{a_1,a_2}$ are products of rising factorials:
\be
   \int\limits_0^\infty x^n \: d\mu_{a_1,a_2}(x)
   \;=\;
   {\Gamma(a_1 + n) \, \Gamma(a_2 + n)  \over \Gamma(a_1) \, \Gamma(a_2)}
   \;=\;
   a_1^{\overline{n}} \, a_2^{\overline{n}}
 \label{def.besselK.moments}
\ee
where $a^{\overline{n}} \eqdef a(a+1) \cdots (a+n-1)$
(see e.g.\ \cite[p.~388]{Watson_44}).

Now fix $a_1, a_2 > 0$ and consider the pair of measures
$(\mu_1, \mu_2) = (\mu_{a_1,a_2}, \mu_{a_1 +1,a_2})$.
Let $P_\bfn(x)$ be the (monic) multiple orthogonal polynomials of type~II
associated to the pair $(\mu_1,\mu_2)$,
and let $\Ptilde_n(x)$ be those polynomials on the stepline:
\be
   \Ptilde_{2k}(x) \;=\; P_{k,k}(x) \;,\qquad
   \Ptilde_{2k+1}(x) \;=\; P_{k+1,k}(x)
   \;.
\ee
Then {Van Assche} and Yakubovich \cite[Theorem~4]{VanAssche_00}
showed that these polynomials satisfy the four-term recurrence\footnote{
   The translation from our notation to theirs is
   $\nu = a_1 - a_2$, $\alpha = a_2 - 1$,
   $b_n = \pi_{n,n}$, $c_n = \pi_{n,n-1}$, $d_n = \pi_{n,n-2}$.
}
\be
   x \, \Ptilde_n(x)
   \;=\;
   \Ptilde_{n+1}(x)
   \:+\: \pi_{n,n} \Ptilde_n(x)
   \:+\: \pi_{n,n-1} \Ptilde_{n-1}(x)
   \:+\: \pi_{n,n-2} \Ptilde_{n-2}(x)
\ee
where
\begin{subeqnarray}
   \pi_{n,n}    & = &  a_1 a_2 \,+\, (2a_1 + 2a_2 - 1)n \,+\, 3 n^2 \\[2mm]
   \pi_{n,n-1}  & = &  n (a_1 +n-1) (a_2 +n-1)(a_1 + a_2 + 3n - 2)   \\[2mm]
   \pi_{n,n-2}  & = &  n (n-1) (a_1 +n-1) (a_1 +n-2) (a_2 +n-1) (a_2 +n-2)
 \label{eq.besselK.pi}
\end{subeqnarray}

On the other hand, P\'etr\'eolle, Zhu and~I have found
\cite[Section~13]{latpath_SRTR},
for all integers $m \ge 1$,
an $m$-branched S-fraction for the ratio of contiguous hypergeometric
series $\FHyper{m+1}{0}$ \cite[Theorem~13.1]{latpath_SRTR}:
namely, if we define the polynomials $P_n^{(m)}(a_1,\ldots,a_m;a_{m+1})$ by
\be
   \sum_{n=0}^\infty P_n^{(m)}(a_1,\ldots,a_m;a_{m+1}) \: t^n
   \;\:=\;\:
   {\FHYPERbottomzero{m+1}{a_1,\ldots,a_{m+1}}{t}
    \over
    \FHYPERbottomzero{m+1}{a_1,\ldots,a_m,a_{m+1}-1}{t}
   }
   \;\,,
 \label{eq.thm.rF0}
\ee
then $P_n^{(m)}(a_1,\ldots,a_m;a_{m+1}) = S_n^{(m)}(\balpha)$
where $S_n^{(m)}$ is the $m$-Stieltjes--Rogers polynomial
and the coefficients $\balpha = (\alpha_i)_{i \ge m}$ are given by
\be
   \balpha
   \;=\;
   a_1 \cdots a_m, \,
   a_2 \cdots a_{m+1}, \,
   a_3 \cdots a_{m+1} (a_1 + 1), \,
   a_4 \cdots a_{m+1} (a_1 + 1)(a_2 + 1), \,
   \ldots
   \;\,.
 \label{eq.thm.rF0.alphas}
\ee
Note that these $\balpha$ can be interpreted
as the products of $m$ successive ``pre-alphas'':
\be
   \balphapre
   \;=\;
   a_1,\ldots,a_{m+1}, a_1+1,\ldots,a_{m+1}+1, a_1+2,\ldots,a_{m+1}+2, \ldots
   \;\,.
 \label{eq.thm.rF0.prealphas}
\ee

In particular, if $a_{m+1} = 1$, then the denominator series $\FHyper{m+1}{0}$
on the right-hand side of \reff{eq.thm.rF0} becomes simply the constant 1,
so that $P_n^{(m)}(a_1,\ldots,a_m;1)$ is simply a product of rising factorials:
\be
   P_n^{(m)}(a_1,\ldots,a_m;1)
   \;=\;
   \prod_{i=1}^m a_i^{\overline{n}}
\ee
(this special case is \cite[Corollary~13.3]{latpath_SRTR}).
Specializing further to $m=2$, we obtain from \reff{eq.thm.rF0.alphas}
\begin{subeqnarray}
   \alpha_{3k+2}  & = &   (a_1 + k) (a_2 + k) \\[1mm]
   \alpha_{3k+3}  & = &   (a_2 + k) (1 + k)  \\[1mm]
   \alpha_{3k+4}  & = &   (1 + k) (a_1 + k + 1)
 \label{eq.3F0.alphas}
\end{subeqnarray}
Then the corresponding production matrix \reff{eq.prop.contraction}
\cite[Propositions~7.2 and 8.2 and eqn.~(7.8)]{latpath_SRTR}
is quadridiagonal with $\pi_{n,n+1} = 1$ and
\begin{subeqnarray}
   \pi_{n,n}    & = &  \alpha_{3n} \,+\, \alpha_{3n+1} \,+\, \alpha_{3n+2}
      \\[1mm]
   \pi_{n,n-1}  & = &  \alpha_{3n-2} \alpha_{3n} \,+\,
                       \alpha_{3n-1} \alpha_{3n} \,+\,
                       \alpha_{3n-1} \alpha_{3n+1}
      \\[1mm]
   \pi_{n,n-2}  & = &  \alpha_{3n-4} \alpha_{3n-2} \alpha_{3n}
 \label{eq.3F0.pi}
\end{subeqnarray}
provided that we make the convention $\alpha_0 = \alpha_1 = 0$.\footnote{
   If $\alpha_{3n}$ and $\alpha_{3n+1}$ are given by polynomial expressions
   in $n$ that do {\em not}\/ vanish when $n=0$,
   then $\pi_{n,n}$ (resp.~$\pi_{n,n-1}$) is given by the corresponding
   polynomial expression {\em plus}\/ a correction term
   proportional to $\delta_{n,0}$ (resp.~$\delta_{n,1}$).
}
The formulae \reff{eq.3F0.alphas} satisfy this convention,
and substituting \reff{eq.3F0.alphas} into \reff{eq.3F0.pi}
gives precisely \reff{eq.besselK.pi}.

In fact, we can go farther and compute the full output matrix $\scro(\Pi)$,
i.e.\ compute the generalized 2-Stieltjes--Rogers polynomials
$S^{(2)}_{n,k}(\balpha)$ for the coefficients $\balpha$
given by \reff{eq.3F0.alphas}.
This was not done in \cite{latpath_SRTR}, but we can do it here:

\begin{proposition}[Generalized 2-Stieltjes--Rogers polynomials
   associated to the rising-factorial moments]
 \label{prop.besselK}
The output matrix $\scro(\Pi) = \sfS^{(2)}(\balpha)$ corresponding to the
production matrix \reff{eq.besselK.pi} is
\be
   S^{(2)}_{n,k}(\balpha)
   \;=\;
   \binom{n}{k} \, (a_1 + k)^{\overline{n-k}} \, (a_2 + k)^{\overline{n-k}}
   \;.
 \label{eq.prop.besselK}
\ee
\end{proposition}

\proof
Let $a_{nk}$ be the right-hand side of \reff{eq.prop.besselK};
we need to show that it satisfies the recurrence \reff{def.iteration.bis}
when the $\pi_{ik}$ are given by \reff{eq.besselK.pi}.
That is, we need to show that
\begin{eqnarray}
   & &
   \hspace*{-7mm}
   \binom{n}{k} \, (a_1 + k)^{\overline{n-k}} \, (a_2 + k)^{\overline{n-k}}
   \;\;=\;
       \nonumber \\[3mm]
   & &
   \binom{n-1}{k-1} \, (a_1 + k-1)^{\overline{n-k}} \, (a_2 + k-1)^{\overline{n-k}}
       \nonumber \\[3mm]
   & & + \;
   \binom{n-1}{k} \, (a_1 + k)^{\overline{n-k-1}} \, (a_2 + k)^{\overline{n-k-1}}
   \; [a_1 a_2 \,+\, (2a_1 + 2a_2 - 1)k \,+\, 3 k^2]
       \nonumber \\[3mm]
   & & + \;
   \binom{n-1}{k+1} \, (a_1 + k+1)^{\overline{n-k-2}} \, (a_2 + k+1)^{\overline{n-k-2}} \:\times    \nonumber \\
   & &
   \hspace*{2.5cm} [(k+1) (a_1 +k) (a_2 +k)(a_1 + a_2 + 3k +1)]
       \nonumber \\[3mm]
   & & + \;
   \binom{n-1}{k+2} \, (a_1 + k+2)^{\overline{n-k-3}} \, (a_2 + k+2)^{\overline{n-k-3}} \:\times    \nonumber \\
   & &
   \hspace*{2.5cm} [(k+2)(k+1) (a_1 +k+1) (a_1 +k) (a_2 +k+1) (a_2 +k)]
   \;.
   \qquad\qquad
\end{eqnarray}
This is a tedious but straightforward computation:
it is convenient to pull out from the right-hand side a factor
${(n-1)! \over (n-k)! \, (k+2)!} \, (a_1 + k+2)^{\overline{n-k-3}} \, (a_2 + k+2)^{\overline{n-k-3}}$
and then evaluate the remaining polynomial expression.
\qed

\bigskip

{\bf Remarks.}
1.  Our definition $(\mu_1, \mu_2) = (\mu_{a_1,a_2}, \mu_{a_1 +1,a_2})$
is manifestly asymmetric between $a_1$ and $a_2$;
nevertheless, the recurrence \reff{eq.besselK.pi}
and the output matrix \reff{eq.prop.besselK}
are symmetric in $a_1 \leftrightarrow a_2$.
The reason is that if we define $\mu'_2 = \mu_{a_1,a_2 +1}$,
then $\mu'_2$ is a linear combination of $\mu_1$ and $\mu_2$
(and vice versa):
\be
   a_1 \, \mu_{a_1 +1,a_2} \:-\: a_2 \, \mu_{a_1,a_2 +1}
   \;=\;
   (a_1 - a_2) \, \mu_{a_1,a_2}
   \;.
\ee
It follows that (as mentioned in Section~\ref{subsec.MOP})
the pairs $(\mu_1, \mu_2)$ and $(\mu_1, \mu'_2)$
give rise to the {\em same}\/ collection of multiple orthogonal polynomials
$P_\bfn(x)$.

Also, the alphas~\reff{eq.3F0.alphas} are asymmetric in $a_1$ and $a_2$,
but they nevertheless give rise to
a production matrix \reff{eq.3F0.pi}/\reff{eq.besselK.pi}
and output matrix \reff{eq.prop.besselK}
that are symmetric in $a_1 \leftrightarrow a_2$.
This is an instance of the nonuniqueness of branched continued fractions
\cite[Sections~3, 10.1, 12.1, 12.2.1 and 13]{latpath_SRTR}.

2.  Explicit expressions for the stepline polynomials $\Ptilde_n(x)$
can be found \cite{BenCheikh_00,Coussement_01};
they are hypergeometric polynomials $\FHyper{1}{2}$.

3. The work of {Van Assche} and Yakubovich \cite{VanAssche_00}
was subsequently generalized by Kuijlaars and Zhang \cite{Kuijlaars_14}
to general $r \ge 2$:
here the measures $\mu_1,\ldots,\mu_r$
have moments that are products of $r$ rising factorials
\cite[eq.~(1.6)]{Kuijlaars_14},
and their densities are expressed in general in terms of Meijer $G$-functions
\cite[eq.~(1.4)]{Kuijlaars_14}.
The stepline polynomials $\Ptilde_n(x)$ are hypergeometric polynomials
$\FHyper{1}{r}$
\cite[eq.~(3.10) and preceding (3.11)]{Kuijlaars_14}
and can also be written in terms of Meijer $G$-functions
\cite[eq.~(3.11)]{Kuijlaars_14}.
The stepline polynomials $\Ptilde_n(x)$ satisfy
an $(r+2)$-term recurrence relation,
for which the coefficients are computed explicitly in
\cite[Corollary~4.2 and Lemma~4.3]{Kuijlaars_14}.
These coefficients presumably coincide with the production matrix
obtained from the $r$-branched S-fraction
\cite[Corollary~13.3]{latpath_SRTR}
via \cite[Propositions~7.2 and 8.2]{latpath_SRTR},
but I have not explicitly checked this for $r > 2$.

4. More generally, one can consider cases in which the moments of
the measures $\mu_1,\ldots,\mu_r$ are ratios of products of rising factorials,
with $p$ factors in the numerator and $q$ factors in the denominator.
Then the ordinary generating function of the moments of $\mu_1$
is a hypergeometric function $\FHyper{p+1}{q}$ with $a_{p+1} = 1$,
and the recurrence relation for the stepline polynomials
can be compared with the branched continued fractions
in \cite[Theorems~14.3, 14.5 and 14.6]{latpath_SRTR}.
In all these branched continued fractions, $r = \max(p,q)$.
Lima and Loureiro have considered the cases
$(p,q) = (2,1)$ \cite{Lima_20a} and $(p,q) = (2,2)$ \cite{Lima_20b};
and Lima \cite{Lima_22} has very recently considered the case of
general $(p,q)$.
\myendremark

\subsection[Bessel $I_\alpha$ weights $\Rightarrow$ Laguerre moments]{Bessel $\bm{I_\alpha}$ weights $\Rightarrow$ Laguerre moments}
   \label{subsec.BesselI}

For real numbers $\alpha \ge -1$ and $x \ge 0$,
define a positive measure $\mu_{\alpha,x}$ on $[0,\infty)$ by
\be
   d\mu_{\alpha,x}(y)
   \;=\;
   \begin{cases}
     e^{-x} \; \FHYPERtopzero{1}{\alpha+1}{xy}
            \: \displaystyle {1 \over \Gamma(\alpha+1)}
            \: y^\alpha \, e^{-y} \, dy
         & \textrm{for $\alpha > -1$}  \\[6mm]
     x \, e^{-(x+y)} \; \FHYPERtopzero{1}{2}{xy} \, dy
         & \textrm{for $\alpha = -1$}
   \end{cases}
\ee
(Here the weight function for $\alpha = -1$ is the
 limit as $\alpha \to -1$ of the ones for $\alpha > -1$.)
The moments of $\mu_{\alpha,x}$ are
\be
   \int\limits_0^\infty y^n \: d\mu_{\alpha,x}(y)
   \;=\;
   \scrl_n^{(\alpha)}(x)
   \;,
\ee
where $\scrl_n^{(\alpha)}(x)$ is the monic unsigned Laguerre polynomial
\be
   \scrl_n^{(\alpha)}(x)
   \;\eqdef\;
   n! \, L_n^{(\alpha)}(-x)
   \;=\;
   \sum_{k=0}^n \binom{n}{k} \, (n+\alpha)^{\underline{n-k}} \, x^k
 \label{def.scrlna}
\ee
and $\rho^{\underline{n}} \eqdef \rho(\rho-1) \cdots (\rho-n+1)$.
Indeed, this is nothing other than the well-known integral representation
for the Laguerre polynomials \cite[Theorem~5.4]{Szego_75},
\be
   \scrl_n^{(\alpha)}(x)
   \;=\;
   n! \, L_n^{(\alpha)}(-x)
   \;=\;
   e^{-x} x^{-\alpha/2}
   \int\limits_0^\infty y^n \: e^{-y} \, y^{\alpha/2} \,
         I_\alpha(2 \sqrt{xy}) \: dy
   \qquad\hbox{for } \alpha > -1
 \label{eq.laguerre.stieltjes}
\ee
where $I_\alpha$ is the modified Bessel function
of the first kind \cite[p.~77]{Watson_44}
\begin{subeqnarray}
   I_\alpha(z)
   & = &
   \sum_{k=0}^\infty {(z/2)^{\alpha+2k} \over k! \, \Gamma(\alpha+k+1)}
        \\[2mm]
   & = &
   \displaystyle {1 \over \Gamma(\alpha+1)} \: (z/2)^\alpha \:
      \FHYPERtopzero{1}{\alpha+1}{z^2/4}
   \;,
 \label{def.BesselI}
\end{subeqnarray}
together with the corresponding limiting formula when $\alpha \to -1$.

Some years ago, Coussement and {Van Assche} \cite{Coussement_03}
studied the multiple orthogonal polynomials of types~I and II
associated to the pair of measures
$(\mu_1,\mu_2) = (\mu_{\alpha,\xi},\mu_{\alpha+1,\xi})$
where $\alpha > -1$ and $\xi > 0$ are fixed parameters.
(They actually used a slightly different normalization,
 so that their moments are $\xi^n \, \scrl_n^{(\alpha)}(\xi)$
 rather than $\scrl_n^{(\alpha)}(\xi)$:
 see \cite[Lemma~1]{Coussement_03}.
 Their $c$ is our $1/\xi$.)
In particular, Coussement and {Van Assche} \cite{Coussement_03}
computed explicitly
the four-term recurrence relation for the multiple orthogonal polynomials
of type~II along the stepline \cite[Theorem~9]{Coussement_03}.
After translating from their normalization to ours,
this four-term recurrence becomes
\be
   x \, \Ptilde_n(x)
   \;=\;
   \Ptilde_{n+1}(x)
   \:+\: \pi_{n,n} \Ptilde_n(x)
   \:+\: \pi_{n,n-1} \Ptilde_{n-1}(x)
   \:+\: \pi_{n,n-2} \Ptilde_{n-2}(x)
\ee
where
\begin{subeqnarray}
   \pi_{n,n}    & = &   (2n+1+\alpha) \,+\, \xi   \\[1mm]
   \pi_{n,n-1}  & = &   n(n+\alpha) \,+\, 2n\xi \\[1mm]
   \pi_{n,n-2}  & = &   n(n-1) \xi
 \label{eq.besselI.pi}
\end{subeqnarray}
This quadridiagonal production matrix $\Pi$
plays a central role in our forthcoming work \cite{latpath_laguerre}
on the coefficientwise Hankel-total positivity of the Laguerre polynomials.
When $\alpha = -1$ (Lah polynomials) it arises from a 2-branched S-fraction,
as found already in \cite{latpath_lah};
when $\alpha = 0$ (rook polynomials) it arises from a
modified 2-branched S-fraction (see \cite{latpath_laguerre}).

Coussement and {Van Assche} \cite{Coussement_03} also gave
an explicit formula for
the multiple orthogonal polynomials of type~II along the stepline
\cite[Theorem~10 and Corollary~2]{Coussement_03}.
After translating from their notation to ours, it is
\be
   \Ptilde_n(x)
   \;=\;
   (-1)^n \sum_{k=0}^n \binom{n}{k} \, \xi^{n-k} \, \scrl_k^{(\alpha)}(-x)
   \;.
\ee
It is curious that Laguerre polynomials occur here too.

\subsection{Final remarks}

I suspect that the foregoing examples are just the tip of the iceberg,
and that the connection between multiple orthogonal polynomials,
production matrices and branched continued fractions will be fruitful
in both directions.
For instance, using known techniques
(such as vector Pearson equations
 \cite{Douak_95,BenCheikh_00,VanAssche_00,Coussement_01,Coussement_03,%
Lima_20a,Lima_20b})
it may be possible to devise new examples of multiple orthogonal polynomials;
these will then automatically provide a production matrix
for the sequence of moments of $\mu_1$;
this production matrix will in turn automatically arise from
a branched J-fraction \cite[Sections~4--8]{latpath_SRTR},
and in some cases this branched J-fraction may arise 
from contraction of a branched S-fraction
\cite[Propositions~7.2 and 7.6]{latpath_SRTR}.
And conversely, combinatorial or algebraic methods leading to new
production matrices or branched continued fractions
may point the way to new examples of multiple orthogonal polynomials.

\section*{Acknowledgments}

This paper arose out of conversations with Walter Van Assche,
to whom I am extremely grateful.
I am also grateful to the organizers of the
15th International Symposium on Orthogonal Polynomials,
Special Functions and Applications
(Hagenberg, Austria, 22--26 July 2019)
for inviting me to give a talk there;
this allowed me to meet Walter and to discover this unexpected
connection between our respective areas of research.
In particular, it was during that conference that Walter and~I discovered
the first example of this connection,
namely, the one shown in Section~\ref{subsec.BesselK}.

I also wish to thank Alex Dyachenko for helpful conversations
and for drawing my attention to the Ph.D.~thesis of Drake \cite{Drake_06},
and H\'elder Lima for many helpful comments on several drafts of this paper.

This research was supported in part by
the U.K.~Engineering and Physical Sciences Research Council grant EP/N025636/1.

\appendix
\section[$\bm{LU}$ factorization for matrices over a commutative ring: \hfill\break
      Proof of Proposition~\ref{prop.LU}]{$\bm{LU}$ factorization for matrices over a commutative ring:
      Proof of Proposition~\ref{prop.LU}}
   \label{app.LU}

If $A$ is a finite or infinite matrix
over a commutative ring
$R$,
we denote by $A_k$
its $k \times k$ leading principal submatrix,
and by $\Delta_k = \det A_k$ the corresponding leading principal minor,
with the convention $\Delta_0 = 1$.

To prove Proposition~\ref{prop.LU}, we will need the following simple fact:

\begin{lemma}
   \label{lemma.uniqueness}
Let $L = (\ell_{ij})_{i,j=1}^n$ be an $n \times n$ lower-triangular matrix
over a commutative ring $R$;
and assume that none of the diagonal elements $\ell_{ii}$
is zero or a divisor of zero.
Then for each vector $b \in R^n$,
the equation $Lx = b$ has at most one solution $x \in R^n$.
%
\end{lemma}

\proof
$\ell_{11} x_1 = b_1$ has at most one solution $x_1$,
since $\ell_{11}$ is neither zero nor a divisor of zero.
Continuing inductively, we see that
$\ell_{ii} x_i = b_i - \sum_{j=1}^{i-1} \ell_{ij} x_j$
has at most one solution $x_i$.
%
\qed

\proofof{Proposition~\ref{prop.LU}}

(a,b,c)
Let $A = (a_{ij})_{i,j=1}^n$ be an $n \times n$ matrix with entries in $R$,
and suppose that we have a factorization $A = LU$
where $L = (\ell_{ij})_{i,j=1}^n$ is lower-triangular
and $U = (u_{ij})_{i,j=1}^n$ is upper-triangular.
It follows that $A_k = L_k U_k$ for all $k$, and hence that
\be
   \Delta_k \;=\; (\det L_k)(\det U_k)
            \;=\  (\ell_{11} u_{11}) (\ell_{22} u_{22})
                      \,\cdots\, (\ell_{kk} u_{kk})
   \;.
 \label{eq.Deltak}
\ee
This proves (a,b,c).

(d) Now suppose that $A = LU$ where $L$ is \emph{unit}-lower-triangular.
As before we have $A_k = L_k U_k$ for all $k$;
and now $\Delta_k = u_{11} u_{22} \cdots u_{kk}$.
Since by hypothesis none of $\Delta_1,\ldots,\Delta_{n-1}$
is zero or a divisor of zero, we can conclude that none of
$u_{11},\ldots,u_{n-1,n-1}$ is zero or a divisor of zero.

We now prove uniqueness by induction on $n$.
The base case $n=1$ is trivial.
Suppose that the result holds for matrices of size $n-1$;
we wish to prove it for $A \in R^{n \times n}$.
Write $\displaystyle A = \begin{pmatrix}
                              A_{n-1}    &  b  \\
                              c^{\rm T}  &  d
                         \end{pmatrix}$,
$\displaystyle L = 
   \begin{pmatrix}
        L_{n-1}    &  0  \\
        \ell^{\rm T}  &  1
   \end{pmatrix}
$
and
$\displaystyle U = 
   \begin{pmatrix}
        U_{n-1}    &  u  \\
        0          &  u_n
   \end{pmatrix}
$.
Then $A = LU$ says that
\vspace*{-2mm}
\begin{subeqnarray}
   A_{n-1}   & = &  L_{n-1} U_{n-1}  \\[1mm]
   b         & = &  L_{n-1} u        \slabel{eq.LDU.b} \\[1mm]
   c^{\rm T} & = &  \ell^{\rm T} U_{n-1}    \slabel{eq.LDU.c} \\[1mm]
   d         & = &  \ell^{\rm T} u \,+\, u_n    \slabel{eq.LDU.d}
     \label{eq.LDU}
\end{subeqnarray}
By the inductive hypothesis, $A_{n-1}$ has a unique $LU$ factorization
$L_{n-1} U_{n-1}$;
and as previously noted,
none of the diagonal elements of $U_{n-1}$ is zero or a divisor of zero.
Clearly \reff{eq.LDU.b} has the unique solution $u = L_{n-1}^{-1} b$.
Moreover, Lemma~\ref{lemma.uniqueness} implies that \reff{eq.LDU.c}
has at most one solution $\ell^{\rm T}$.
Then $u_n = d - \ell^{\rm T} u$ is determined as well.
This proves~(d).

(e)
The existence proof is also by induction on $n$.
The base case $n=1$ is trivial;
and for the inductive step we again use \reff{eq.LDU}.
By the inductive hypothesis, $A_{n-1}$ has a unique $LU$ factorization
$L_{n-1} U_{n-1}$,
and moreover $U_{n-1}$ is invertible (because $\Delta_{n-1}$ is).
Then (\ref{eq.LDU}b--d) have the unique solution
$u = L_{n-1}^{-1} b$, $\ell^{\rm T} = c^{\rm T} U_{n-1}^{-1}$
and $u_n = d - c^{\rm T} U_{n-1}^{-1} L_{n-1}^{-1} b$.
\qed

\end{document}